\documentclass[final,onefignum,onetabnum]{siamart190516}

\hypersetup{
  pdftitle={Analysis of the multiplicative Schwarz method for matrices with a special block structure},
  pdfauthor={Carlos Echeverr{\'i}a, J{\"o}rg Liesen and Petr Tich{\'y}}
  pdfkeywords={multiplicative Schwarz method, iterative methods, convergence analysis, singularly perturbed problems, Shishkin mesh discretization, block diagonal dominance}
}

\title{Analysis of the multiplicative Schwarz method for matrices with a special
block structure\thanks{The work of C. Echeverr{\'i}a was partially supported by the Berlin Mathematical School and the Einstein Center for Mathematics Berlin. The work of P. Tich{\'y} was supported by the Grant Agency of the Czech Republic under the grants  no. 17-04150J and 20-01074S.
}}

\author{Carlos Echeverr{\'i}a\footnotemark[4], J\"{o}rg Liesen\footnotemark[4] \and Petr Tich{\'y}\footnotemark[2]}

\headers{Analysis of the multiplicative Schwarz method for special block matrices}{C.~Echeverr{\'i}a, J.~Liesen and P. Tich{\'y}}

\usepackage{times}
\usepackage{amsmath}
\usepackage{amssymb}
\usepackage{graphicx}
\usepackage{epstopdf}
\usepackage{color}
\usepackage{pgfplots}
\pgfplotsset{compat=1.10}

\newcommand{\R}{\mathbb{R}}
\newtheorem{example}[theorem]{Example}

\newcommand{\klein}[1]{\scriptscriptstyle #1}
\newcommand{\Tr}{\ensuremath{\klein{T}}}
\newcommand{\entryvE}{\ensuremath{b_{\klein{H}}}}
\newcommand{\entryvW}{\ensuremath{c_{\klein{H}}}}
\newcommand{\entryvC}{\ensuremath{a_{\klein{H}}}}
\newcommand{\entryvN}{\ensuremath{e_{\klein{H}}}}
\newcommand{\entryvS}{\ensuremath{d_{\klein{H}}}}
\newcommand{\entrywE}{\ensuremath{b_{\klein{h}}}}
\newcommand{\entrywW}{\ensuremath{c_{\klein{h}}}}
\newcommand{\entrywC}{\ensuremath{a_{\klein{h}}}}
\newcommand{\entrywN}{\ensuremath{e_{\klein{h}}}}
\newcommand{\entrywS}{\ensuremath{d_{\klein{h}}}}
\newcommand{\entryzE}{\ensuremath{b}}
\newcommand{\entryzW}{\ensuremath{c}}
\newcommand{\entryzC}{\ensuremath{a}}
\newcommand{\entryzN}{\ensuremath{e}}
\newcommand{\entryzS}{\ensuremath{d}}

\newcommand{\matAH}{\ensuremath{A_{\klein{H}}}}
\newcommand{\matAHhat}{\ensuremath{\widehat{A}_{\klein{H}}}}
\newcommand{\matBH}{\ensuremath{B_{\klein{H}}}}
\newcommand{\matCH}{\ensuremath{C_{\klein{H}}}}
\newcommand{\matAhhat}{\ensuremath{\widehat{A}_{\klein{h}}}}
\newcommand{\matAh}{\ensuremath{A_{\klein{h}}}}
\newcommand{\matBh}{\ensuremath{B_{\klein{h}}}}
\newcommand{\matCh}{\ensuremath{C_{\klein{h}}}}
\newcommand{\matA}{\ensuremath{A}}
\newcommand{\matB}{\ensuremath{B}}
\newcommand{\matC}{\ensuremath{C}}

\newcommand{\piOne}{\ensuremath{\Pi^{\klein{(1)}}}}
\newcommand{\piTwo}{\ensuremath{\Pi^{\klein{(2)}}}}
\newcommand{\pOne}{\ensuremath{P^{\klein{(1)}}}}
\newcommand{\pTwo}{\ensuremath{P^{\klein{(2)}}}}

\newcommand{\atopfrac}[2]{\genfrac{}{}{0pt}{}{#1}{#2}}

\begin{document}

\maketitle

\renewcommand{\thefootnote}{\fnsymbol{footnote}}

\footnotetext[4]{Institute of Mathematics, TU Berlin, Stra{\ss}e des 17. Juni 136, 10623 Berlin, Germany.\\
e-mail: {\tt \{echeverria,liesen\}@math.tu-berlin.de}.}
\footnotetext[2]{Faculty of Mathematics and Physics, Charles University, Prague,
Czech Republic. \\
e-mail: {\tt ptichy@karlin.mff.cuni.cz}.}

\begin{abstract}
We analyze the convergence of the (algebraic) multiplicative Schwarz method
applied to linear algebraic systems with matrices having a special block
structure that arises, for example, when a (partial) differential equation is
posed and discretized on a domain that consists of two subdomains with an
overlap. This is a basic situation in the context of domain decomposition
methods. Our analysis is based on the algebraic structure of the Schwarz
iteration matrices, and we derive error bounds that are based on the block
diagonal dominance of the given system matrix. Our analysis does not assume
that the system matrix is symmetric (positive definite), or has the $M$- or
$H$-matrix property. Our approach is motivated by and significantly generalizes
an analysis for a special one-dimensional model problem given
in~\cite{EchLieSzyTic18}.
\end{abstract}



\section{Introduction}\label{sec:intro}
The (algebraic) multiplicative Schwarz method, sometimes also called the
Schwarz alternating method, is a stationary iterative method for solving large
and sparse linear algebraic systems
\begin{equation}\label{eq:linsys}
\mathcal{A}x=b.
\end{equation}
In each step of the method, the current iterate is multiplied by an iteration
matrix that is the product of several factors, where each factor corresponds to an
inversion of only a restricted part of the matrix $\mathcal{A}$. In the context of the
numerical solution of discretized differential equations, the restrictions of
the matrix correspond to different parts of the computational domain. This
motivates the name ``local solve'' given to each factor, which is also used in a purely algebraic
setting.

The convergence theory for the multiplicative Schwarz method is well
established for important matrix classes including symmetric positive definite
matrices and nonsingular $M$-matrices~\cite{BenFroNabSzy01}, symmetric
indefinite~\cite{FroNabSzy08,FroSzy14} and
semidefinite~\cite{NabSzy06} matrices,
and $H$-matrices~\cite{BruPedSzy04}. The derivation of convergence results for
these matrix classes is usually based on splittings of ${\cal A}$. No systematic
convergence theory exists however for general nonsymmetric matrices.

An important source of nonsymmetric linear algebraic systems is the
discretization of singularly perturbed convection-diffusion problems. In a
recent paper we have studied the convergence of the multiplicative Schwarz
method for a one-dimensional model problem in this
context~\cite{EchLieSzyTic18}. The system matrices in this problem are usually
nonsymmetric, nonnormal, ill-conditioned, and in particular not in one of the
classes considered in~\cite{BenFroNabSzy01,BruPedSzy04,FroNabSzy08,FroSzy14}.
Moreover, our derivations are not based on matrix splittings, but on the
off-diagonal decay of the matrix inverses, which in turn is implied by diagonal
dominance. From a broader point of view our results show why a convergence
theory for the multiplicative Schwarz method for ``general'' matrices will most
likely remain elusive: Even in the simple model problem considered
in~\cite{EchLieSzyTic18}, the convergence of the method strongly depends on the
problem parameters and on the chosen discretization, and while the method
rapidly converges in some cases, it even diverges in others.

Ideally, the study of well chosen model problems involves a
``generalizability aspect'' in the sense that the obtained results shall give
insight into more general problems; see for example the interesting discussion
about the nature of model problems in~\cite{KraParSte83}. This is also true for
the approach in~\cite{EchLieSzyTic18}, since it motivated our analysis in this
paper, which may be considered a significant generalization of the previous
work. Here we analyze the convergence of the multiplicative Schwarz method for
linear algebraic systems with matrices of the form
\begin{equation}\label{eq:blockmat}
\mathcal{A}=\left[
  \begin{array}{ccc}
             \matAHhat       & e_{m}\otimes \matBH   &    0            \\
    e_{m}^{\Tr}\otimes \matC  &   A      & e_{1}^{\Tr}\otimes \matB \\
                   0         & e_{1}\otimes\matCh &        \matAhhat  \\
  \end{array}
\right]\;\in\;\R^{N(2m+1)\times N(2m+1)},
\end{equation}
with $\matAHhat, \matAhhat \in \R^{Nm\times Nm}$,
$\matA, \matB, \matC, \matBH, \matCh \in \R^{N\times N}$, and
the canonical basis vectors $e_1,e_m\in\R^{m}$. We will usually think of
$\matAHhat, \matAhhat \in \R^{Nm\times Nm}$ as matrices consisting of $m$
blocks of size $N\times N$. After deriving general expressions for the norms of
the multiplicative Schwarz iteration matrices for systems of the form
\eqref{eq:linsys}--\eqref{eq:blockmat}, we derive actual error bounds only for
the case when the blocks $\matAHhat$ and $\matAhhat$ of ${\cal A}$ are block
tridiagonal.

Such matrices arise naturally when a differential equation is posed and
discretized inside a domain $\Omega$ that is divided into two subdomains
$\Omega_1$ and $\Omega_2$ with an overlap. In this context the first $m+1$
block rows in the matrix $\mathcal{A}$ correspond to the unknowns in the domain
$\Omega_1$, the last $m+1$ block rows correspond to the unknowns in the domain
$\Omega_2$, and the middle block row corresponds to the unknowns in the
overlap. The underlying assumption here is that in each of the two subdomains
we have the same number of unknowns. This assumption is made for simplicity of
our exposition. Extensions to other block sizes are possible, but would require
even more technicalities; see our discussion in \Cref{sec:conclusions}.
We point out that the model problems studied in~\cite{EchLieSzyTic18} are of the
form \eqref{eq:linsys}--\eqref{eq:blockmat} with $N=1$. While (usual) diagonal
dominance of tridiagonal matrices is one of the main tools
in~\cite{EchLieSzyTic18}, the derivation of error bounds here relies on recent
results on block diagonal dominance of block tridiagonal matrices
from~\cite{EchLieNab18}.

The paper is organized as follows. In \Cref{sec:schwarz} we state the
multiplicative Schwarz method for linear algebraic systems of the form
\eqref{eq:linsys}--\eqref{eq:blockmat}, and in \Cref{sec:structure} we
study the algebraic structure and the norm of its iteration matrices. In
\Cref{sec:bounds} we derive error bounds for the method when the matrix
${\cal A}$ is block tridiagonal and block diagonally dominant. We apply these
error bounds in \Cref{sec:application} to a two-dimensional discretized
convection-diffusion model problem. Finally, in \Cref{sec:conclusions}
we summarize the main results of the paper and briefly discuss possible
generalizations and alternative applications of our approach.

\section{The multiplicative Schwarz method}\label{sec:schwarz}

The multiplicative Schwarz method for solving linear algebraic systems of the
form \eqref{eq:linsys}--\eqref{eq:blockmat} can naturally be based on two local
solves using the top and the bottom $N(m+1)\times N(m+1)$ block of ${\cal A}$,
respectively. More precisely, the restriction operators of the method are
\[
R_1 \equiv \left[ I_{N(m+1)} \quad 0 \right]
\quad \mbox{and} \quad
R_2 \equiv \left[ 0 \quad I_{N(m+1)} \right],
\]
which are both of size $N(m+1)\times N(2m+1)$.
The corresponding restrictions of $\mathcal{A}$ are
\[
A_1 \equiv R_1\mathcal{A}R_1^{\Tr}=
\left[
 \begin{array}{cc}
             \matAHhat      &  e_{m} \otimes \matBH \\
     e_{m}^{\Tr} \otimes C  &             A
 \end{array}
\right] \;\,\mbox{and}\;\,
A_2 \equiv R_2\mathcal{A}R_2^{\Tr}=
\left[
 \begin{array}{cc}
             A          &  e_{1}^{\Tr} \otimes B \\
     e_{1} \otimes C_h  &       \matAhhat
 \end{array}
\right],
\]
which are both of size $N(m+1)\times N(m+1)$. We now define the two projections
\begin{equation}\label{eq:Pi}
P_i\equiv R_i^{\Tr}A_i^{-1}R_i\mathcal{A}\;\in\;
\R^{N(2m+1)\times N(2m+1)},
\quad i=1,2.
\end{equation}
Then their complementary projections
\[
Q_i\equiv I-P_i\;\in\;\R^{N(2m+1)\times N(2m+1)},\quad i=1,2,
\]
yield the multiplicative Schwarz iteration matrices
\begin{equation}\label{eqn:Tij}
T_{12}\equiv Q_2Q_1\quad\mbox{and}\quad T_{21}\equiv Q_1Q_2.
\end{equation}
If in the context of discretized differential equations the top and bottom
blocks of ${\cal A}$ correspond to the unknowns in the domains $\Omega_1$ and
$\Omega_2$, respectively, then a multiplication with the matrix $T_{ij}$
corresponds to first performing a local solve on $\Omega_i$, and then on
$\Omega_j$.

We now fix either $T_{12}$ or $T_{21}$, and choose an initial vector
$x^{(0)}\in \R^{N(2m+ 1)}$. Then the (algebraic, one-level) multiplicative
Schwarz method is defined by
\begin{equation}\label{eq:schwarz}
x^{(k+1)}=T_{ij} x^{(k)}+v, \quad k=0,1,2,\dots\,.
\end{equation}
The vector $v\in \R^{N(2m+ 1)}$ is defined to make the method
consistent. The consistency condition is given by $x=T_{ij}x+v$, and this
yields
\[
v = (I-T_{ij})x = (P_1+P_2-P_jP_i)x,
\]
which is (easily) computable since
\[
P_ix=R_i^{\Tr}A_i^{-1}R_i\mathcal{A}x=R_i^{\Tr}A_i^{-1}R_ib,\quad i=1,2.
\]
The error of the multiplicative Schwarz iteration~\eqref{eq:schwarz} is given by
\begin{equation}\label{eq:schwarz_err}
e^{(k+1)}=x-x^{(k+1)}=(T_{ij}x + v)- (T_{ij} x^{(k)}+v)=T_{ij} e^{(k)},
\quad k=0,1,2,\dots,
\end{equation}
and hence $e^{(k+1)}=T_{ij}^{k+1} e^{(0)}$ by induction.
For any consistent norm $\|\cdot\|$, we therefore have the error bound
\begin{equation}\label{eq:error}
\|e^{(k+1)}\|\leq \|T_{ij}^{k+1}\|\,\|e^{(0)}\|.
\end{equation}
In the following we will use the structure of the iteration matrices $T_{ij}$
to derive bounds on the norms $\|T_{ij}^{k+1}\|$. This will lead to stronger
convergence results than the standard approach in the convergence analysis
of stationary iterative methods, which uses submultiplicativity, i.e.,
$\|T_{ij}^{k+1}\|\leq \|T_{ij}\|^{k+1}$, and then bounds $\|T_{ij}\|$.
As mentioned above, our approach is motivated by~\cite{EchLieSzyTic18}, but
since we consider matrices with block rather than scalar entries, neither
the analysis nor the results from~\cite{EchLieSzyTic18} are directly
applicable here.

\section{Structure and norms of the iteration matrices} \label{sec:structure}

Let us have a closer look at the structure of the matrices $T_{ij}$. A direct
computation based on \eqref{eq:Pi} shows that
\[
P_1 =
\left[ \begin{array}{c}
I_{N(m+1)} \\
    0
\end{array} \right]
\matA_1^{-1}
\left[ \begin{array}{c|c|c}
\matA_1 & e_{m+1}\otimes B & 0
\end{array} \right]
=
\left[\begin{array}{ccc}
I_{N(m+1)} & \matA_1^{-1}(e_{m+1}\otimes B) & 0\\
    0 &    0                               & 0
\end{array}\right],
\]
and
\[
P_2 =
\left[\begin{array}{c}
    0    \\
I_{N (m+1)}
\end{array}\right]
\matA_2^{-1}
\left[\begin{array}{c|c|c}
 0 &  e_{1}\otimes C & \matA_2
\end{array}\right]
=
\left[\begin{array}{ccc}
0 &            0                        &  0      \\
0 & \matA_2^{-1}(e_{1}\otimes C) & I_{N (m+1)}
\end{array}\right],\nonumber
\]
where $e_1,e_{m+1}\in\R^{m+1}$. We see that both $P_1$ and $P_2$ have
exactly $N(m+1)$ linearly independent columns, and hence
\[
\mathrm{rank}(P_1) = \mathrm{rank}(P_2)=N(m+1).
\]
Moreover, the complementary projections are
\[
\hspace*{-0.5em}
Q_1  =
\left[
  \begin{array}{ccc}
    0 & -\matA_1^{-1}(e_{m+1}\otimes B)  &  \\
     & I_{N}                            &  \\
     &                                 & I_{N(m-1)}
  \end{array}
\right],\;
Q_2  =
\left[
  \begin{array}{ccc}
    I_{N(m-1)} &                             &  \\
              & I_{N}                         &  \\
              & -\matA_2^{-1}(e_{1}\otimes C) & 0
  \end{array}
\right],
\]
and we have
\[
\mathrm{rank}(Q_1) = \mathrm{rank}(Q_2)=Nm.
\]
In order to simplify the notation we write
\begin{equation}\label{eq:pi_and_p}
\left[
  \begin{array}{c}
    \pOne \\ \piOne
  \end{array}
\right]
\equiv
\matA_1^{-1}(e_{m+1}\otimes B)
\quad\mathrm{and}\quad
\left[
  \begin{array}{c}
    \piTwo \\ \pTwo
  \end{array}
\right]
\equiv
\matA_2^{-1}(e_{1}\otimes C),
\end{equation}
where $\Pi^{\klein{(i)}} \in \mathbb{R}^{N\times N}$, and
\[
P^{\klein{(i)}} =
\left[\left(P_1^{\klein{(i)}}\right)^{\Tr},\ldots,
\left(P_{m}^{\klein{(i)}}\right)^{\Tr}\right]^{\Tr}
\in \R^{Nm\times N}\quad\mathrm{with}\quad P^{\klein{(i)}}_j
\in \R^{N\times N},\; j=1,\dots,m.
\]
Then
\[
Q_1  =
\left[\begin{array}{cccc}
             0 &     & -\pOne     &             \\
                 & 0 & -\piOne        &             \\
                 &     &      I_N       &             \\
                 &     &                & I_{N (m-1)}
      \end{array}\right]\quad\mbox{and}\quad
Q_2 =
\left[\begin{array}{cccc}
I_{N(m-1)}  &                &     &                 \\
            & I_N            &     &                 \\
            & -\piTwo        & 0 &                 \\
            & -\pTwo         &     & 0
\end{array}\right],
\]
so that
\begin{eqnarray} \label{eqn:T12}
T_{12} &=& Q_2Q_1 =
\left[
  \begin{array}{ccc}
    0  &      -\pOne               &  0     \\
    0  & \piTwo \pOne_{m}          &  0     \\
    0  & \pTwo \pOne_{m} &  0
  \end{array}
\right] \\
 &=&\left[
  \begin{array}{c}
    -\pOne                  \\
	  \piTwo\pOne_{m}   \\
	  \pTwo \pOne_{m}  \\
  \end{array}
\right]
\left[
  \begin{array}{c|c|c}
    0_{N(m+1)} & I_N & 0_{N(m-1)}
  \end{array}
\right]
\equiv V_1(e_{m+2}^{\Tr}\otimes I_N), \nonumber
\end{eqnarray}
and
\begin{eqnarray} \label{eqn:T21}
T_{21}&=&Q_1Q_2 =
\left[
  \begin{array}{ccc}
    0  &  \pOne \pTwo_{1}  &  0     \\
    0  & \piOne\pTwo_{1}           &  0     \\
    0  &  -\pTwo                   &  0
  \end{array}
\right] \\
&=&\left[
  \begin{array}{c}
    \pOne \pTwo_{1}    \\
	  \piOne\pTwo_{1}           \\
	  -\pTwo                    \\
  \end{array}
\right]
\left[
  \begin{array}{c|c|c}
    0_{N(m-1)} & I_N & 0_{N(m+1)}
  \end{array}
\right]
\equiv V_2(e_{m}^{\Tr}\otimes I_N), \nonumber
\end{eqnarray}
where $e_m,e_{m+2}\in\R^{2m+1}$ and $V_1,V_2\in\R^{N(2m+1)\times N}$.

Using these representations of the matrices $T_{ij}$, we can obtain the
following generalization of~\cite[Proposition~4.1 and
Corollary~4.2]{EchLieSzyTic18}.

\begin{lemma}\label{lem:powers}
In the notation established above we have $\mathrm{rank}(T_{ij})\leq N$ and
\begin{equation}\label{eq:Tkp1}
T_{12}^{k+1} =  V_1 \left( \pTwo_{1}\pOne_{m} \right)^{k}
(e_{m+2}^{\Tr}\otimes I_N),
\quad
T_{21}^{k+1} =  V_2\left(\pOne_{m}\pTwo_{1}\right)^{k}(e_{m}^{\Tr}\otimes I_N),
\end{equation}
for all $k\geq 0$.
\end{lemma}
\begin{proof}
We only consider the matrix $T_{12}$; the proof for $T_{21}$ is analogous.
The result about the rank is obvious from~\eqref{eqn:T12}. We denote
$E_{m+2}\equiv e_{m+2}^{\Tr}\otimes I_N$, then $T_{12}=V_1E_{m+2}$, and it is
easy to see that
\[
T_{12}^{k+1} =V_1\left(E_{m+2}V_1\right)^{k}E_{m+2},\quad
\mbox{for all $k\geq 0$.}
\]
Now
\[
E_{m+2}V_1=(e_{m+2}^{\Tr}\otimes I_N)
\left[
\begin{array}{c}
  -\pOne                         \\
   \piTwo \pOne_{m}              \\
   \pTwo_{1} \pOne_{m}           \\
   \pTwo_{2:m} \pOne_{m}
\end{array}
\right]=\pTwo_{1}\pOne_{m},
\]
which shows the first equality in \eqref{eq:Tkp1}.
\end{proof}

The next result generalizes~\cite[Lemma~4.3]{EchLieSzyTic18} and gives
expressions for some of the block entries of the matrices $T_{ij}$, which will
be essential in our derivations of error bounds in the following section.

\begin{lemma} \label{lem:pp}
Suppose that the matrices $\matAHhat, \matAhhat \in \R^{Nm\times Nm}$ in
\eqref{eq:blockmat} are nonsingular, and denote
$\matAHhat^{-1}=[Z^{\klein{(H)}}_{ij}]$ and
$\matAhhat^{-1}=[Z^{\klein{(h)}}_{ij}]$ with
$Z^{\klein{(H)}}_{ij},Z^{\klein{(h)}}_{ij}\in \R^{N\times N}$.
Then, in the notation established above,
\[
\left[
\begin{array}{c}
\pOne \\ \piOne
\end{array}
\right] =
\left[
\begin{array}{c}
-Z^{\klein{(H)}}_{1:m,m}\matBH \\  I_{N}
\end{array}
\right] \piOne,
\quad
\piOne = \left( \matA - \matC Z^{\klein{(H)}}_{mm}
\matBH \right)^{-1}\matB, 
\]
and
\[
\left[
\begin{array}{c}
\piTwo\\
\pTwo
\end{array}
\right] =
\left[
\begin{array}{c}
I_{N} \\
-Z^{\klein{(h)}}_{1:m,1}\matCh
\end{array}
\right] \piTwo,
\quad
\piTwo =\left( \matA - \matB Z^{\klein{(h)}}_{11} \matCh \right)^{-1}\matC.
\]
\end{lemma}

\begin{proof}
From \eqref{eq:pi_and_p} we know that $\pOne$, $\pTwo$, $\piOne$, and $\piTwo$
solve the linear algebraic systems
\[
\hspace*{-1em}
\left[
\begin{array}{cc}
        \matAHhat         & e_{m} \otimes \matBH \\
e_{m}^{\Tr} \otimes \matC & \matA
\end{array}
\right]
\left[
\begin{array}{c}
\pOne \\
\piOne
\end{array}
\right]=
\left[
\begin{array}{c}
    0 \\
 \matB
\end{array}
\right],
\quad
\left[
\begin{array}{cc}
        \matA          & e_{1}^{\Tr} \otimes \matB \\
 e_{1} \otimes \matCh  &  \matAhhat
\end{array}
\right]
\left[
\begin{array}{c}
\piTwo \\
\pTwo
\end{array}
\right]
= \left[
\begin{array}{c}
\matC \\
 0
\end{array}
\right].\nonumber
\]
Hence the expressions for  $\pOne$, $\pTwo$, $\piOne$, and $\piTwo$ can be
obtained using Schur complements; see, e.g.,~\cite[\S~0.7.3]{HornJohn12}.
\end{proof}

In order to bound norms of powers of the
iteration matrices $T_{ij}$ we have to decide first
which matrix norm should be taken. In the following we use a general induced
matrix norm $\|\cdot\|$ which can be considered for square as well as for
rectangular matrices. Note that an induced matrix norm for square matrices
is submultiplicative and satisfies $\|I\|=1$.

\begin{lemma}\label{lem:genbounds}
In the notation established above, for any induced matrix norm we have
\begin{equation}\label{eqn:powers}
\|T_{ij}^{k+1}\| \leq \rho_{ij}^{k} \|T_{ij}\|,\quad\mbox{for all $k\geq 0$,}
\end{equation}
where
\begin{equation}\label{eqn:rho}
\rho_{12}=\|Z_{11}^{\klein{(h)}}\matCh \piTwo
Z_{mm}^{\klein{(H)}}\matBH\piOne\|
\quad\mbox{and}\quad
\rho_{21} =\|Z_{mm}^{\klein{(H)}}\matBH \piOne
Z_{11}^{\klein{(h)}}\matCh \piTwo\|.
\end{equation}
\end{lemma}

\begin{proof}
We only consider the matrix $T_{12}$; the proof for $T_{21}$ is analogous.
Taking norms in \eqref{eq:Tkp1} yields
\[
\|T_{12}^{k+1}\|= \|V_1(\pTwo_{1}\pOne_{m})^{k}E_{m+2}\|
\leq \rho_{12}^k \|V_1\| \|E_{m+2}\|
\quad\mbox{with}\quad \rho_{12} \equiv  \|\pTwo_{1}\pOne_{m}\|,
\]
and where
$E_{m+2}$ is defined as in the proof of \Cref{lem:powers}.
Note that $\|E_{m+2}\|=\|I_N\|=1$ and
\[
\|T_{12}\|=\max_{\|x\|=1}\|T_{12}x\|=\max_{\|x\|=1}\|V_1E_{m+2}x\|
=\max_{\|y\|=1}\|V_1y\|=\|V_1\|,
\]
which yields the bound on $\|T_{12}^{k+1}\|$ in \eqref{eqn:powers}.\pagebreak

Finally, the equality
$\|\pTwo_{1}\pOne_{m}\|=\|Z_{11}^{\klein{(h)}}\matCh \piTwo
Z_{mm}^{\klein{(H)}}\matBH\piOne\|$ in \eqref{eqn:rho} follows directly
from \Cref{lem:pp}.
\end{proof}

So far our analysis considered general nonsingular blocks $\matAHhat$ and
$\matAhhat$ in the matrix ${\cal A}$ in \eqref{eq:blockmat}, and combining
\eqref{eq:error} and \eqref{eqn:powers} gives a general error bound for the
multiplicative Schwarz method in terms of certain blocks of ${\cal A}$ and the
inverses of $\matAHhat$ and $\matAhhat$. Note that using the
submultiplicativity of the matrix norm $\|\cdot\|$, which at this point is
still a general induced norm, both convergence factors $\rho_{12}$ and
$\rho_{21}$ can be bounded by
\begin{equation}\label{eqn:rhobound}
\rho_{ij}\leq \|Z_{11}^{\klein{(h)}}\matCh\|\,
\|Z_{mm}^{\klein{(H)}}\matBH\|\, \| \piOne\|\,  \| \piTwo\|.
\end{equation}
In order to derive a quantitative error bound from the terms on the
right hand side, we have to make additional assumptions on $\matAHhat$ and
$\matAhhat$. One possible choice of such assumptions is considered in the next
section.

\section{Error bounds for the block tridiagonal case}\label{sec:bounds}

We are most interested in the analysis of the multiplicative Schwarz method for
linear algebraic systems that arise in certain discretizations of partial
differential equations, and we will therefore consider
\begin{equation}\label{eqn:tridiag}
\matAHhat=\mathrm{tridiag}(\matCH,\matAH,\matBH)
\quad\text{and}\quad
\matAhhat=\mathrm{tridiag}(\matCh,\matAh,\matBh).
\end{equation}
Additionally, we will assume that the matrices
\begin{equation}\label{eqn:nonsing}
\matAH,\matBH,\matCH,\matA,\matB,\matC,\matAh,\matBh,\matCh \in
\R^{N\times N}
\;\;\mbox{are nonsingular,}
\end{equation}
and that the matrix ${\cal A}$ is {\em row block diagonally dominant} in the
sense of \cite[Definition~2.1]{EchLieNab18}, i.e., that
\begin{eqnarray}
 \|\matAH^{-1}\matBH\|+\|\matAH^{-1}\matCH\| &\leq& 1,\nonumber\\
 \|A^{-1}B\|+\|A^{-1}C\| &\leq& 1,\label{eqn:dominant}\\
 \|\matAh^{-1}\matBh\|+\|\matAh^{-1}\matCh\| &\leq& 1.\nonumber
\end{eqnarray}
Note that because of \eqref{eqn:nonsing}, each of the norms on the left hand
sides of these inequalities is \emph{strictly} less than one.

Both $\matAH$ and $\matAh$ satisfy all assumptions
of~\cite[Theorem~2.6]{EchLieNab18}. A minor modification of the first equation
in the proof of that theorem (namely multiplying both sides of this equation by
$\matCh$ or $\matBH$ before taking norms) shows that the blocks of the
inverses, i.e., $\matAHhat^{-1}=[Z^{\klein{(H)}}_{ij}]$ and
$\matAhhat^{-1}=[Z^{\klein{(h)}}_{ij}]$, satisfy
\begin{equation}\label{eqn:decay}
\|Z^{\klein{(h)}}_{i1}\matCh\|\leq \|Z^{\klein{(h)}}_{11}\matCh\|
\quad \mbox{and} \quad
\|Z^{\klein{(H)}}_{im}\matBH\|\leq \|Z^{\klein{(H)}}_{mm}\matBH\|,
\quad i=1,\dots,m.
\end{equation}
Moreover, as shown in the proof of that theorem, the equations
\begin{eqnarray*}
Z_{11}^{\klein{(h)}} &=& (\matAh-\matBh M_h)^{-1} =
(I-\matAh^{-1}\matBh M_h)^{-1}\matAh^{-1},\\
Z_{mm}^{\klein{(H)}} &=& (\matAH-\matCH L_H)^{-1}=
(I-\matAH^{-1}\matCH L_H)^{-1}\matAH^{-1}
\end{eqnarray*}
hold for some matrices $M_h,L_H\in\R^{N\times N}$ with $\|M_h\|\leq 1$
and $\|L_H\|\leq 1$; see~\cite[equation (2.20)]{EchLieNab18}. The precise
definition of $M_h$ and $L_H$ is not important here.

The four matrices that appear on the right hand side of \eqref{eqn:rhobound}
are now given by
\begin{eqnarray*}
Z_{11}^{\klein{(h)}}\matCh &=&
(I-\matAh^{-1}\matBh M_h)^{-1}\matAh^{-1}\matCh,\\
Z_{mm}^{\klein{(H)}}\matBH &=&
(I-\matAH^{-1}\matCH L_H)^{-1}\matAH^{-1}\matBH,\\
\piTwo &=& (I-\matA^{-1}\matB Z_{11}^{\klein{(h)}}\matCh)^{-1}\matA^{-1}\matC,\\
\piOne &=& (I-\matA^{-1}\matC Z_{mm}^{\klein{(H)}}\matBH)^{-1}\matA^{-1}\matB.
\end{eqnarray*}
Since $\|\matAh^{-1}\matBh M_h\|\leq \|\matAh^{-1}\matBh\| <1$, we can use the
Neumann series to obtain
\begin{eqnarray*}
\|(I-\matAh^{-1}\matBh M_h)^{-1}\| &=&
\left\|\sum_{k=0}^\infty (\matAh^{-1}\matBh M_h)^k\right\|\\
&\leq& \sum_{k=0}^\infty \|\matAh^{-1}\matBh\|^k =
\frac{1}{1-\|\matAh^{-1}\matBh\|}.
\end{eqnarray*}
Similarly, $\|\matAH^{-1}\matCH L_H\|\leq \|\matAH^{-1}\matCH\| <1$ implies that
\[
\|(I-\matAH^{-1}\matCH L_H)^{-1}\|\leq \frac{1}{1-\|\matAH^{-1}\matCH\|},
\]
and hence
\begin{eqnarray}
\|Z_{11}^{\klein{(h)}}\matCh\| &\leq&
\frac{\|\matAh^{-1}\matCh\|}{1-\|\matAh^{-1}\matBh\|}
\equiv \eta_h \leq 1, \label{eqn:etah} \\
\|Z_{mm}^{\klein{(H)}}\matBH\| &\leq&
\frac{\|\matAH^{-1}\matBH\|}{1-\|\matAH^{-1}\matCH\|}
\equiv \eta_H \leq 1.\label{eqn:etaH}
\end{eqnarray}
Using \eqref{eqn:etah} and \eqref{eqn:etaH} yields
\[
\|\matA^{-1}\matB Z_{11}^{\klein{(h)}}\matCh\|\leq \eta_h
\|\matA^{-1}\matB\| <1\quad\mbox{and}\quad
\|\matA^{-1}\matC Z_{mm}^{\klein{(H)}}\matBH\|\leq \eta_H
\|\matA^{-1}\matC\| <1,
\]
and another application of the Neumann series shows that
\begin{equation}\label{eqn:Sh}
\|\piTwo\| \leq
\frac{\|\matA^{-1}\matC\|}{1-\eta_h \|\matA^{-1}\matB\|}\leq 1
\quad \mbox{and} \quad
\|\piOne\| \leq
\frac{\|\matA^{-1}\matB\|}{1-\eta_H \|\matA^{-1}\matC\|}\leq 1.
\end{equation}
In summary, we have the following result.
\begin{lemma}\label{lem:rho}
In the notation established above, the convergence factors of the multiplicative
Schwarz method satisfy
\begin{equation}\label{eqn:rhoij}
\rho_{ij}\leq
\frac{\eta_h \|\matA^{-1}\matC\|}{1-\eta_h \|\matA^{-1}\matB\|}\;
\frac{\eta_H\|\matA^{-1}\matB\|}{1-\eta_H\|\matA^{-1}\matC\|},
\end{equation}
where each of the factors on the right hand side is less than or equal to one.
\end{lemma}

Let us illustrate the bound from \Cref{lem:rho} on a simple example.

\begin{example}{\rm
Let $m\geq 1$ be given, $N=2m+1$, and consider the matrix
\[
{\cal A}\equiv {\rm tridiag}(-I,W,-I)\in{\mathbb R}^{N(2m+1)\times N(2m+1)},
\]
where
$I\in {\mathbb R}^{N\times N}$ and
\[
    W \equiv {\rm tridiag}(-1,4,-1)= 4\left(I - \frac14 S\right)\in {\mathbb R}^{N\times N}, \quad
    S\equiv {\rm tridiag}(1,0,1) \in {\mathbb R}^{N\times N}.
\]
It is well known that ${\cal A}$ is the result
of a standard finite difference discretization of the 2D Poisson equation on
the unit square and with Dirichlet boundary conditions. In our notation,
${\cal A}$ is of the form \eqref{eq:blockmat} and \eqref{eqn:tridiag} with
\[
\matAHhat=\matAhhat=\mathrm{tridiag}(-I,W,-I) \in{\mathbb R}^{Nm\times Nm},
\]
$\matBH=B=\matBh=\matCH=C=\matCh=-I$, $\matAH=A=\matAh=W$, and
\[
H = h = \frac{1}{N+1}.
\]
\pagebreak

\noindent Let us first investigate row block diagonal dominance of $\mathcal{A}$
with respect to the matrix $2$-norm and $\infty$-norm.
It both cases it holds that $\| S \| \leq 2$, since $\| S \|_2 \leq \| S \|_\infty = 2$.
By expanding $W^{-1}$ using the Neumann series
we get the bound
\begin{eqnarray}\label{eq:boundW}
\| W^{-1} \| &=&
\frac14 \left\| \sum_{k=0}^{\infty} \left(\frac14 S\right)^k \right\|
\leq
 \frac14 \sum_{k=0}^{\infty} \left(\frac{\|S\|}{4}\right)^k
 = \frac{1}{4-\|S\|}\,,
\end{eqnarray}
and, therefore, $\| W^{-1} \| \leq \frac12$.

For the $\infty$-norm the inequality in \eqref{eq:boundW} is strict since
$\|S^k\|_\infty < \|S\|_\infty^k$ for $k>N/2$, and hence
\[
   \|W^{-1}\|_2 \leq \|W^{-1}\|_\infty < \frac{1}{4-\|S\|_\infty} = \frac12.
\]
As a consequence, ${\cal A}$ is strictly row block diagonally dominant with
respect to both considered matrix norms; see the conditions~\eqref{eqn:dominant}. Note that
${\cal A}$ is only weakly row diagonally dominant in the classical (scalar)
sense.

Let us now concentrate on bounding the convergence factor $\rho_{ij}$
of the multiplicative Schwarz method; see \eqref{eqn:rhoij}. Using the definitions \eqref{eqn:etah} and \eqref{eqn:etaH} we obtain
\[
    \eta_h=\eta_H=\frac{\| W^{-1} \| }{1-\| W^{-1} \|}
\]
and \eqref{eqn:rhoij} yields
\[
\rho_{ij}\leq
\left(\frac{\frac{\| W^{-1} \| }{1-\| W^{-1} \|} \|W^{-1}\|}{1-\frac{\| W^{-1} \| }{1-\| W^{-1} \|} \|W^{-1}\|}\right)^2.
\]
Since this bound is an increasing function of
$\| W^{-1} \|$ in $(0,1)$, we can use \eqref{eq:boundW} to obtain
\begin{equation}\label{eq:boundS}
\rho_{ij}\leq
\left(\frac{1  }{11-7\|S\| + \| S\|^2}\right)^2.
\end{equation}
For the $\infty$-norm, the inequality in \eqref{eq:boundS}
is strict and the right-hand side
is exactly one. For the $2$-norm, the eigenvalues of $S$
are known explicitly, it holds that
\[
\| S\|_2 = 2-2\sin^2\left(\frac{\pi h}{2}\right)
\]
and, therefore,
\[
\rho_{ij}\leq\left(\frac{1}{1 + 6\sin^2\left(\frac{\pi h}{2}\right) +
4\sin^4\left(\frac{\pi h}{2}\right)}\right)^2
\approx
\frac{1}{1+3\pi^2h^2}
\]
for small values of $h$.
Thus, the convergence factor of the multiplicative Schwarz method for
both considered norms is less than one, regardless of the choice
of $m$. But note that for $N\rightarrow \infty$, and hence $h\rightarrow 0$,
we can in both cases expect that $\rho_{ij}\rightarrow 1$.
}\end{example}

\medskip
In order to bound the error norm of the multiplicative Schwarz method, see
\eqref{eq:error}, \eqref{eqn:powers}, and \eqref{eqn:rhoij}, it remains to
bound $\|T_{ij}\|$. Let us first note that because of the equivalence of
matrix norms, there exists a constant $c$ such that
\[
   \|T_{ij}\| \leq c\, \|T_{ij}\|_\infty,
\]
where $c$ can depend on the size of $T_{ij}$.\pagebreak

Now we bound $\|T_{ij}\|_{\infty}$.
From \eqref{eqn:T12} and \eqref{eqn:T21} we see that
\begin{eqnarray}
\|T_{12}\|_{\infty} &=&\max\{ \|\pOne\|_{\infty},\,
\|\piTwo\pOne_{m}\|_{\infty},\,
\|\pTwo \pOne_{m}\|_{\infty}\},\label{eqn:T12a}\\
\|T_{21}\|_{\infty} &=&\max\{ \|\pTwo\|_{\infty},\, \|\piOne\pTwo_1\|_{\infty},\,
\|\pOne\pTwo_1\|_{\infty}\},\label{eqn:T21a}
\end{eqnarray}
and \Cref{lem:pp} yields
\begin{align*}
\pOne &= -Z_{1:m,m}^{\klein{(H)}}\matBH \piOne,
&
\pTwo &= -Z_{1:m,1}^{\klein{(h)}}\matCh \piTwo,
\\
\piTwo\pOne_{m} &=-\piTwo Z_{mm}^{\klein{(H)}}\matBH \piOne,
&
\piOne\pTwo_1 &=-\piOne Z_{11}^{\klein{(h)}}\matCh \piTwo,
\\
\pTwo \pOne_{m} &= Z_{1:m,1}^{\klein{(h)}} \matCh \piTwo Z_{mm}^{\klein{(H)}}
\matBH \piOne,
&
\pOne \pTwo_1 &= Z_{1:m,m}^{\klein{(H)}}\matBH \piOne Z_{11}^{\klein{(h)}}
\matCh\piTwo.
\end{align*}
Using \eqref{eqn:decay} we can bound the $\infty$-norms of these matrices as
follows:
\begin{eqnarray*}
\|\pOne\|_{\infty} &=& \max\{ \|Z_{1m}^{\klein{(H)}}\matBH\piOne\|_{\infty},
\dots, \|Z_{mm}^{\klein{(H)}}\matBH\piOne\|_{\infty} \}\\
&\leq& \|Z_{mm}^{\klein{(H)}}\matBH\|_{\infty} \| \piOne \|_{\infty},\\
\|\piTwo\pOne_{m}\|_{\infty} &=&
\|\piTwo Z_{mm}^{\klein{(H)}}\matBH\piOne\|_{\infty}\\
&\leq& \|Z_{mm}^{\klein{(H)}}\matBH\|_{\infty} \|\piOne\|_{\infty}
\|\piTwo\|_{\infty},\\
\|\pTwo \pOne_{m}\|_{\infty} &=& \max\{ \|Z_{11}^{\klein{(h)}}\matCh \piTwo
Z_{mm}^{\klein{(H)}}\matBH \piOne\|_{\infty},
\dots,\|Z_{m1}^{\klein{(h)}}\matCh\piTwo Z_{mm}^{\klein{(H)}}\matBH
\piOne\|_{\infty} \}\\
&\leq& \|Z_{mm}^{\klein{(H)}}\matBH\|_{\infty} \|\piOne\|_{\infty}
\|Z_{11}^{\klein{(h)}}\matCh\|_{\infty}   \|\piTwo\|_{\infty},
\end{eqnarray*}
and
\begin{eqnarray*}
\|\pTwo\|_{\infty} &=&\max\{ \|Z_{11}^{\klein{(h)}}\matCh\piTwo\|_{\infty},
\dots, \|Z_{m1}^{\klein{(h)}}\matCh\piTwo\|_{\infty}\}\\
&\leq& \|Z_{11}^{\klein{(h)}}\matCh\|_{\infty} \| \piTwo \|_{\infty},\\
\|\piOne\pTwo_1\| & =&
\| \piOne Z_{11}^{\klein{(h)}}\matCh\piTwo \|_{\infty} \\
&\leq&  \|Z_{11}^{\klein{(h)}}\matCh\|_{\infty} \| \piTwo \|_{\infty}
\|\piOne\|_{\infty}, \\
\|\pOne \pTwo_1\|_{\infty} &=& \max\{ \|Z_{1m}^{\klein{(H)}}\matBH \piOne
Z_{11}^{\klein{(h)}}\matCh \piTwo\|_{\infty},
\dots,\| Z_{mm}^{\klein{(H)}}\matBH \piOne
Z_{11}^{\klein{(h)}}\matCh\piTwo\|_{\infty} \}\\
&\leq& \|Z_{mm}^{\klein{(H)}}\matBH\|_{\infty} \|\piOne\|_{\infty}
\|Z_{11}^{\klein{(h)}}\matCh\|_{\infty}   \|\piTwo\|_{\infty}.
\end{eqnarray*}
The individual terms on the right hand sides of previous inequalities are all
less or equal than one. Therefore, the maximum of the first three bounds is
$\|Z_{mm}^{\klein{(H)}}\matBH\|_{\infty} \| \piOne \|_{\infty}$, and the
maximum of the second three bounds is $\|Z_{11}^{\klein{(h)}}\matCh\|_{\infty}
\| \piTwo \|_{\infty}$. Hence, using \eqref{eqn:T12a} and \eqref{eqn:T21a}, and
\eqref{eqn:etah}, \eqref{eqn:etaH}, \eqref{eqn:Sh} we obtain
\[
	\| T_{12} \|_{\infty} \leq \|Z_{mm}^{\klein{(H)}}\matBH\|_{\infty}
  \| \piOne \|_{\infty} \leq
	\frac{\eta_{H,\klein{\infty}}\|\matA^{-1}\matB\|_{\infty}}
  {1-\eta_{H,\klein{\infty}} \|\matA^{-1}\matC\|_{\infty}}
\]
and
\[
	\| T_{21} \|_{\infty} \leq \|Z_{11}^{\klein{(h)}}\matCh\|_{\infty}
  \| \piTwo \|_{\infty} \leq
	\frac{\eta_{h,\klein{\infty}}\|\matA^{-1}\matC\|_{\infty}}
  {1-\eta_{h,\klein{\infty}} \|\matA^{-1}\matB\|_{\infty}},
\]
where $\eta_{h,\klein{\infty}}$ and $\eta_{H,\klein{\infty}}$ are defined as in
\eqref{eqn:etah} and \eqref{eqn:etaH} using the $\infty$-norm,
\begin{equation}\label{eqn:etas}
\eta_{h,\klein{\infty}} \equiv \frac{\|\matAh^{-1}\matCh\|_{\infty}}
{1-\|\matAh^{-1}\matBh\|_{\infty}},
\qquad
\eta_{H,\klein{\infty}} \equiv \frac{\|\matAH^{-1}\matBH\|_{\infty}}
{1-\|\matAH^{-1}\matCH\|_{\infty}}.
\end{equation}
Combining these bounds with \Cref{lem:genbounds} and \Cref{lem:rho}
gives the following convergence result.

\begin{theorem}\label{thm:main}
Suppose that ${\cal A}$ as in \eqref{eq:blockmat} has blocks as in
\eqref{eqn:tridiag} that satisfy \eqref{eqn:nonsing}--\eqref{eqn:dominant}.
Then the errors of the  multiplicative Schwarz method \eqref{eq:schwarz}
applied to the linear algebraic system \eqref{eq:linsys} satisfy
\[
\frac{\|e^{(k+1)}\|}{\|e^{(0)}\|} \leq
\left(\frac{\eta_h \|\matA^{-1}\matC\|}{1-\eta_h \|\matA^{-1}\matB\|}\;
\frac{\eta_H\|\matA^{-1}\matB\|}
{1-\eta_H\|\matA^{-1}\matC\|}\right)^k\|T_{ij}\|,
\quad k=0,1,2,\dots,
\]
where $\eta_h$ and $\eta_H$ are defined in \eqref{eqn:etah}--\eqref{eqn:etaH}.
Moreover,
\[
	\| T_{12} \| \leq
	c\, \frac{\eta_{H,\klein{\infty}}\|\matA^{-1}\matB\|_{\infty}}
  {1-\eta_{H,\klein{\infty}} \|\matA^{-1}\matC\|_{\infty}}  ,\qquad
	\| T_{21} \| \leq
	c\, \frac{\eta_{h,\klein{\infty}}\|\matA^{-1}\matC\|_{\infty}}
  {1-\eta_{h,\klein{\infty}} \|\matA^{-1}\matB\|_{\infty}},
\]
where $\eta_h^{\klein{(\infty)}}$ and $\eta_H^{\klein{(\infty)}}$ are given by
\eqref{eqn:etas}, and $c$ is a constant such that $\|T_{ij}\| \leq c\,
\|T_{ij}\|_\infty$.
\end{theorem}

In many practical applications the matrices $\matAHhat$ and $\matAhhat$ of the
form \eqref{eqn:tridiag} are not only row block diagonally dominant,
see \eqref{eqn:dominant}, but also {\em column block diagonally dominant},
i.e., they satisfy the conditions
\begin{equation}\label{eqn:dominant2}
 \|\matBH\matAH^{-1}\|+\|\matCH\matAH^{-1}\| \leq 1,\quad
 \|\matBh\matAh^{-1}\|+\|\matCh\matAh^{-1}\| \leq 1;
\end{equation}
see \Cref{def:coldom} in Appendix~A. Analogously to the row
block diagonally dominant case described above, a proof of
\Cref{thm:coldom} along the lines of the proof of~\cite[Theorem~2.6]{EchLieNab18}
shows that if $\matAHhat$ and $\matAhhat$
satisfy the conditions \eqref{eqn:dominant2}, then
\begin{eqnarray*}
Z_{11}^{\klein{(h)}} &=& (\matAh-\tilde{L}_h \matCh)^{-1} =
\matAh^{-1}(I-\tilde{L}_h\matCh \matAh^{-1})^{-1},\\
Z_{mm}^{\klein{(H)}} &=& (\matAH-\tilde{M}_H \matBH )^{-1}=
\matAH^{-1}(I-\tilde{M}_H\matBH \matAH^{-1})^{-1},
\end{eqnarray*}
for some matrices $\tilde{L}_h, \tilde{M}_H\in\R^{N\times N}$ with
$\|\tilde{L}_h\|\leq 1$ and $\|\tilde{M}_H\|\leq 1$. Analogously, using the
Neumann series we obtain the bounds
\[
\|Z_{11}^{(h)}\|	\leq	\frac{\|\matAh^{-1}\|}{1-\|\matCh\matAh^{-1}\|},\qquad
\|Z_{mm}^{(H)}\|	\leq	\frac{\|\matAH^{-1}\|}{1-\|\matBH\matAH^{-1}\|},
\]
so that
\[
\|Z_{11}^{(h)}\matCh\|\leq\frac{\|\matAh^{-1}\|\|\matCh\|}
{1-\|\matCh\matAh^{-1}\|},\qquad
\|Z_{mm}^{(H)}\matBH\|	\leq	\frac{\|\matAH^{-1}\| \|\matBH\|}
{1-\|\matBH\matAH^{-1}\|}.
\]
Therefore, if $\matAHhat$ and $\matAhhat$ satisfy the conditions
\eqref{eqn:dominant} \emph{and} \eqref{eqn:dominant2}, then
\begin{equation}\label{eqn:Z11min}
\|Z_{11}^{(h)}\matCh\|\leq \min \left\{ \frac{\|\matAh^{-1}\matCh\|}
{1-\|\matAh^{-1}\matBh\|},  \frac{\|\matAh^{-1}\|\|\matCh\|}
{1-\|\matCh\matAh^{-1}\|} \right\} \equiv \eta_h^{\klein{\min}},
\end{equation}
and
\begin{equation}\label{eqn:Zmmmin}
\|Z_{mm}^{(H)}\matBH\|	\leq
\min \left\{\frac{\|\matAH^{-1}\matBH\|}{1-\|\matAH^{-1}\matCH\|},
\frac{\|\matAH^{-1}\| \|\matBH\|}{1-\|\matBH\matAH^{-1}\|}\right\} \equiv
\eta_H^{\klein{\min}}.
\end{equation}
The value of $\eta_H^{\klein{\min}}$ can be much smaller than $\eta_H$ for
example if $\|\matBH\|\ll \|\matCH\|$.
Since we only improved bounds \eqref{eqn:etah} and \eqref{eqn:etaH}
on $\|Z_{11}^{(h)}\matCh\|$ and $\|Z_{mm}^{(H)}\matBH\|$,
we can formulate a version of \Cref{thm:main} where we just replace
$\eta_h$ and $\eta_H$ with $\eta_h^{\klein{\min}}$ and $\eta_H^{\klein{\min}}$.

\begin{theorem}\label{thm:main2}
Suppose that ${\cal A}$ as in \eqref{eq:blockmat} has blocks as in
\eqref{eqn:tridiag} that satisfy \eqref{eqn:nonsing}, \eqref{eqn:dominant}, and
\eqref{eqn:dominant2}. Then the errors of the  multiplicative Schwarz method
\eqref{eq:schwarz} applied to the linear algebraic system \eqref{eq:linsys}
satisfy
\[
\frac{\|e^{(k+1)}\|}{\|e^{(0)}\|} \leq
\left(\frac{\eta_h^{\klein{\min}} \|\matA^{-1}\matC\|}{1-\eta_h^{\klein{\min}}
\|\matA^{-1}\matB\|}\;
\frac{\eta_H^{\klein{\min}}\|\matA^{-1}\matB\|}{1-
\eta_H^{\klein{\min}}\|\matA^{-1}\matC\|}\right)^k\|T_{ij}\|,
\quad k=0,1,2,\dots,
\]
where $\eta_h^{\klein{\min}}$ and $\eta_H^{\klein{\min}}$ are defined in
\eqref{eqn:Z11min} and \eqref{eqn:Zmmmin}.\pagebreak

\noindent Moreover,
\[
	\| T_{12} \| \leq
	c\, \frac{\eta_{H,\klein{\infty}}^{\klein{\min}}
  \|\matA^{-1}\matB\|_{\infty}}{1-\eta_{H,\klein{\infty}}^{\klein{\min}}
  \|\matA^{-1}\matC\|_{\infty}}  ,\qquad
	\| T_{21} \| \leq
	c\, \frac{\eta_{h,\klein{\infty}}^{\klein{\min}}
  \|\matA^{-1}\matC\|_{\infty}}{1-\eta_{h,\klein{\infty}}^{\klein{\min}}
  \|\matA^{-1}\matB\|_{\infty}},
\]
where $\eta_{h,\klein{\infty}}^{\klein{\min}}$ and
$\eta_{H,\klein{\infty}}^{\klein{\min}}$ are given by \eqref{eqn:Z11min} and
\eqref{eqn:Zmmmin} using the $\infty$-norm, and $c$ is a constant such that
$\|T_{ij}\| \leq c\, \|T_{ij}\|_\infty$.
\end{theorem}

In the next section we will apply these results in the study of a model problem
that can be considered a two-dimensional generalization of the one-dimensional
problem studied in~\cite{EchLieSzyTic18}.

\section{Application to a discretized convection-diffusion equation}\label{sec:application}

We consider the following convection-diffusion model problem with Dirichlet
boundary conditions,
\begin{equation}\label{eq:bvp2}
-\epsilon  \Delta u
+ u_y +\beta u = f  \text{ in } \Omega=(0,1)\times(0,1),
\quad \quad
u=g\,\text{  on  }\, \partial\Omega.
\end{equation}
Here the scalar-valued function $u(x,y)$ represents the concentration of a
transported quantity, $\epsilon >0$ is the (scalar) diffusion parameter, and
$\beta \geq 0$ is the (scalar) reaction parameter. We assume that the problem
is \emph{convection dominated}, i.e., that $\epsilon\ll 1$.
Moreover, we assume that the parameters of the problem are chosen so that the
solution $u(x,y)$ has one boundary layer at $y=1$.
A common approach for discretizing such a problem is to use Shishkin meshes.
This technique has been described in detail, for example, in the
articles~\cite[Section~5]{Sty05} and~\cite{KopOri10}, as well as in the
book~\cite{MilOriShi96}. We therefore only briefly summarize the facts that are
relevant for our analysis.

\begin{figure}
\hspace*{-2em}
\begin{minipage}[b]{0.5\textwidth}
\begin{center}
\begin{tikzpicture}[scale = 0.4]
\draw (0,0) rectangle (10, 10);
\node[left] at (0,10) {$1$};
\draw (0,8) node[left] {$1-\tau_y$} -- (10,8);
\node[left] at (0,0) {$0$};
\node[below] at (0,0) {$0$};
\node[below] at (10,0) {$1$};
\draw[black,fill=black] (0,8) circle (0.5ex);
\node at (5,4) {$\Omega_1$};
\node at (5,9) {$\Omega_2$};
\end{tikzpicture}
\end{center}
\end{minipage}
\hspace*{1em}
\begin{minipage}[b]{0.5\textwidth}
\begin{center}
\begin{tikzpicture}[scale = 0.4]
\draw (0,0) rectangle (10, 10);

\foreach \point in {2,4,6,8,8.5,9,9.5}
\draw (0,\point) -- (10,\point);
\node[left] at (0,10) {$y_M$};
\node[left] at (0,8)  {$y_{M/2}$};
\node[left] at (0,0)  {$y_0$};

\foreach \point in {2,4,6,8}
\draw (\point,0) -- (\point,10);
\node[below] at (0,0) {$x_0$};
\node[below] at (10,0) {$x_N$};

\draw[black,fill=black] (0,8) circle (0.5ex);
\draw (5,-1.25) node[above] {$H_x$};
\draw (10,2.5) node[right] {$H_y$};
\draw (10,8.7) node[right] {$h_y$};
\end{tikzpicture}
\end{center}
\end{minipage}
\caption{Division of the domain and Shishkin mesh for the problem
\eqref{eq:bvp2} with one boundary layer at $y=1$.}
\label{fig:2Ddomain2}
\end{figure}
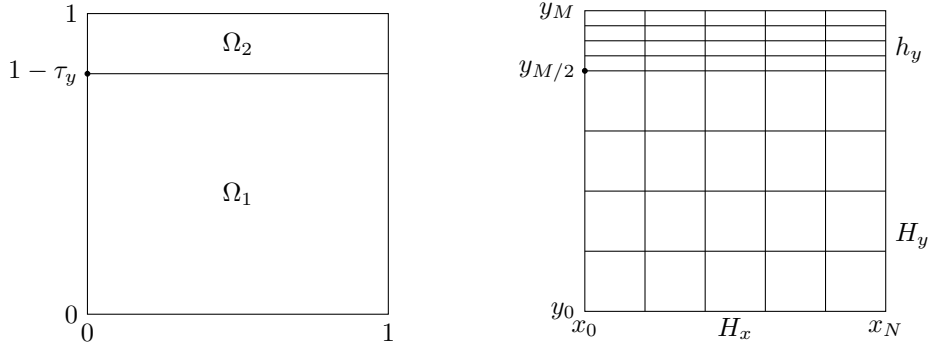

The main idea of the Shishkin mesh is to use piecewise equidistant meshes in
each part of the domain, and to resolve the boundary layer by using a finer
mesh close to $y=1$. In our example the \emph{transition point} at which the
mesh changes from coarse to fine in the $y$-direction is defined by
\begin{equation}\label{eq:tau}
1-\tau_y,\quad\mbox{where}\quad \tau_y \equiv \min\left\{\frac12,2\epsilon  \ln M\right\},
\end{equation}
and where the even integer $M\geq 4$ denotes the number of mesh intervals in
the $y$-direction. The assumption $\epsilon\ll 1$ implies that $\tau_y\ll 1$,
so that $1-\tau_y$ will be very close to $y=1$, unless~$M$ is chosen extremely
large. The domain $\overline{\Omega}$ is thus decomposed into the two
overlapping subdomains
\[
\Omega_{1}=[0,1]\times[0,1-\tau_y]\quad\mbox{and}\quad
\Omega_{2}=[0,1]\times[1-\tau_y,1];
\]
see the left part of \Cref{fig:2Ddomain2}.

Let the integer $N\geq 3$ denote the number of equidistant intervals used in
the $x$-direction; see the right part of \Cref{fig:2Ddomain2}.
In the literature on Shishkin mesh discretizations (including the references cited
above) we usually find $N=M$, but being able to choose $N\neq M$ gives
some additional flexibility in the construction of the mesh.

We denote the mesh width in the $x$-direction by $H_x$, and the widths before
and after the transition point in the $y$-direction by $H_y$ and $h_y$, i.e.,
\begin{equation}\label{eq:nHh}
H_x\equiv \frac{1}{N},\quad H_y\equiv \frac{2(1-\tau_y)}{M},\quad h_y\equiv \frac{2\tau_y}{M}.
\end{equation}
Then the nodes of the Shishkin mesh are given by
\begin{eqnarray}
& \left\{(iH_x,y_j)\,:\,i=0,\dots,N; j=0,\dots, M \right\}\,
\subset\,\overline{\Omega},\quad\mbox{where}\label{eq:meshnodes}\\
& y_j\equiv
\begin{cases} jH_y &\mbox{for}\;\; j=0,\ldots,M/2, \\
1-(M-j)h_y &\mbox{for}\;\; j=M/2+1,\ldots,M. \end{cases}\nonumber
\end{eqnarray}
The ratio between the different mesh sizes in the $y$-direction is given by
\[\frac{h_y}{H_y}=\frac{\tau_y}{1-\tau_y}=\tau_y+{\cal O}(\tau_y^2) \ll 1.\]

Using the standard upwind finite difference discretization of \eqref{eq:bvp2}
and the lexicographical line ordering of the unknowns yields a linear algebraic
system with $\mathcal{A}$ as in \eqref{eq:blockmat} and the submatrices
$\matAHhat$  and $\matAhhat$ having the block tridiagonal structure
\eqref{eqn:tridiag}. Each block row of $\mathcal{A}$ corresponds to one row of
unknowns in the mesh, and hence each block is of size \linebreak$(N-1)\times (N-1)$;
cf. the right part of \Cref{fig:2Ddomain2}. Note that due to the
standard notation in the context of Shishkin mesh discretizations we here slightly
depart from the notation in \eqref{eq:blockmat} and \eqref{eqn:tridiag},
where the individual blocks are of size $N\times N$. In the $y$-direction
there are $M-1$ interior nodes, and one of them is at the transition point
$y_{M/2}=1-\tau_y$. In the notation of \eqref{eq:blockmat}, we have $m=M/2-1$,
and $\mathcal{A}$ is of size $(N-1)(M-1)\times (N-1)(M-1)$.

Following the description in~\cite{LinSty01,Sty05}, we see that the (nonzero)
off-diagonal blocks of ${\cal A}$ are given by
\[
C_H = d_H I,\quad C = d I,\quad C_h = d_h I,\quad
B_H = e_H I,\quad B = e I,\quad B_h = e_h I,
\]
where
\begin{align}
	&\entryvS  =  - \frac{\epsilon}{H_y^{2}}-\frac{1}{H_y},&
	&\entryzS  =  - \frac{2\epsilon}{H_y(H_y+h_y)}-\frac{1}{H_y},&
	&\entrywS  =  - \frac{\epsilon}{h_y^{2}}-\frac{1}{h_y},\label{eq:upwind}\\
    &\entryvN  =  - \frac{\epsilon}{H_y^{2}},&
    &\entryzN  =  - \frac{2\epsilon}{h_y(H_y+h_y)},&
    &\entrywN  =  - \frac{\epsilon}{h_y^{2}}.\nonumber
\end{align}
and $A_H=\text{tridiag}(\entryvW,\entryvC,\entryvE)$, ${A}=\text{tridiag}(\entryzW,\entryzC,\entryzE)$, $A_h=\text{tridiag}(\entrywW,\entrywC,\entrywE)$, where
\begin{align}
    &\entryvW  =  -\frac{\epsilon}{H_x^{2}},& 
    &\entryvC  =  \frac{2\epsilon}{H_x^{2}} +\frac{2\epsilon}{H_y^{2}}
    +\frac{1}{H_y}+\beta,
    &\entryvE  =  -\frac{\epsilon}{H_x^{2}},\nonumber  \\
    &\entryzW  =  -\frac{\epsilon}{H_x^{2}} ,&
    &\entryzC  =  \frac{2\epsilon}{H_x^{2}} +\frac{2\epsilon}{H_yh_y}
    +\frac{1}{H_y}+\beta,
    &\entryzE  =  -\frac{\epsilon}{H_x^{2}}, \label{eq:upwind2}\\
    &\entrywW  =  -\frac{\epsilon}{H_x^{2}},& 
    &\entrywC  =  \frac{2\epsilon}{H_x^{2}} +\frac{2\epsilon}{h_y^{2}}
    +\frac{1}{h_y}+\beta,
    &\entrywE  =  -\frac{\epsilon}{H_x^{2}}.\nonumber
\end{align}
We will now show that for this model problem the assumptions of
\Cref{thm:main2} are satisfied.

\begin{lemma}\label{lem:BDDcon-dif}
All nonzero blocks of the matrix $\mathcal{A}$ described above are nonsingular.
Moreover, for the matrix $\infty$-norm the matrix $\mathcal{A}$ satisfies the
conditions \eqref{eqn:dominant}, i.e., it is row block diagonally dominant, and
the submatrices $\matAHhat$ and $\matAhhat$ satisfy the conditions
\eqref{eqn:dominant2}, i.e., they are column block diagonally dominant.
\end{lemma}

\begin{proof}
Note that all (nonzero) off-diagonal entries of ${\cal A}$ are negative,
and that the diagonal entries $a_H$, $a$, $a_h$ are positive. Moreover,
\[
a_H + b_H + c_H + d_H + e_H
= a + b + c + d + e
= a_h + b_h + c_h + d_h + e_h = \beta \geq 0
\]
It is thus easy to see that all nonzero blocks of ${\cal A}$ are nonsingular.

To prove \eqref{eqn:dominant} and \eqref{eqn:dominant2} for the $\infty$-norm,
we just need to show that
\begin{equation}\label{eqn:shishkinblocks}
 |e_H+d_H| \|\matAH^{-1}\|_\infty \leq 1,\quad
|e+d| \|A^{-1}\|_\infty \leq 1,\quad
|e_h+d_h| \|\matAh^{-1}\|_\infty \leq 1,
\end{equation}
and hence we need to bound the $\infty$-norms of matrices
$\matAH^{-1}$, $\matA^{-1}$, and $\matAh^{-1}$.

First note that for any nonsingular matrix $\mathcal{M}$ and an
induced matrix norm we have
\[
\|\mathcal{M}^{-1}\|=\max_{\|v\|=1}\left\Vert \mathcal{M}^{-1}\left(\frac{\mathcal{M}v}{\|\mathcal{M}v\|}\right)\right\Vert =\frac{1}{\min_{\|v\|=1}\|\mathcal{M} v\|}.
\]
Therefore, if $\|\mathcal{M} v\|\geq\gamma>0$ for any unit norm vector $v$,
then $\|\mathcal{M}^{-1}\|\leq{\gamma^{-1}}$.

Second, suppose that $\mathcal{M}$ is a strictly diagonally dominant
tridiagonal Toeplitz matrix
$\mathcal{M}=\mathrm{tridiag}(\hat{c},\hat{a},\hat{b})$,
where $\hat{a}>0$, $\hat{b}<0$, $\hat{c}<0$, and
$
\hat{a}+\hat{b}+\hat{c}>0.
$
We would like to bound $\|\mathcal{M} v\|_{\infty}$ for any unit norm vector $v$
from below. If $\|v\|_{\infty}=1$, then there is an index $i$ such that
$|v_{i}|=1$. Without loss of generality we can assume that $v_{i}=1$, because
changing the sign of the vector does not change $\|\mathcal{M} v\|_{\infty}$.
Defining $v_{0}=0$ and $v_{n+1}=0$ we obtain
\[
\|\mathcal{M} v\|_{\infty}\geq|v_{i-1}\hat{c}+\hat{a}+v_{i+1}\hat{b}|\geq \hat{a}+\hat{b}+\hat{c},
\]
and therefore
\begin{equation}\label{eqn:boundT}
\|\mathcal{M}^{-1}\|_{\infty}\leq\frac{1}{\hat{a}+\hat{b}+\hat{c}}.
\end{equation}

In order to prove \eqref{eqn:shishkinblocks}, we now apply the bound
\eqref{eqn:boundT} to matrices $\matAH$, $\matA$, and $\matAh$, which are
strictly diagonally dominant tridiagonal Toeplitz matrices with the required
sign pattern. For $\matAH$ we get
\[
    |e_H + d_H| \|\matAH^{-1}\|_\infty \leq \frac{|e_H + d_H|}{a_H+b_H+c_H}
    = \frac{|e_H + d_H|}{|e_H + d_H| + \beta} \leq 1,
\]
and the other inequalities in \eqref{eqn:shishkinblocks} follow analogously.
\end{proof}

\Cref{lem:BDDcon-dif} ensures that the assumptions of \Cref{thm:main2}
are satisfied for matrices arising from the discretization of the problem
\eqref{eq:bvp2} that we have described above (namely, upwind differences
on a Shishkin mesh). We therefore obtain the following convergence
result for the multiplicative Schwarz method.

\begin{corollary}\label{cor:conv-diff}
Consider the linear algebraic system \eqref{eq:linsys} arising from the upwind discretization of the boundary value problem \eqref{eq:bvp2}
on the Shishkin mesh \eqref{eq:meshnodes}. Then the errors of the
multiplicative Schwarz method \eqref{eq:schwarz} applied to this system
satisfy
\begin{equation}\label{eq:cor:conv-diff_1}
  \frac{\|e^{(k+1)}\|_{\infty}}{\|e^{(0)}\|_{\infty}} \leq
  \rho^k\|T_{ij}\|_{\infty},
  \quad k=0,1,2,\dots,
\end{equation}
where
\begin{equation}\label{eq:cor:conv-diff_2}
\rho\equiv \frac{\epsilon}{\epsilon+H_y},\qquad
\| T_{12} \|_{\infty} \leq \rho
  \quad\text{and}\quad \| T_{21} \|_{\infty} \leq  1.\pagebreak
\end{equation}
\end{corollary}

\begin{proof}
To prove the bounds \eqref{eq:cor:conv-diff_1}--\eqref{eq:cor:conv-diff_2}
we apply \Cref{thm:main2} with the $\infty$-norm and bound the
factors $\eta_{h,\infty}^{\klein{\min}}$ and $\eta_{H,\infty}^{\klein{\min}}$
that correspond to the discretization scheme~\eqref{eq:upwind}.

Using $|d_{h}|>|e_{h}|$ we get
\[
  \eta_{h,\infty}^{\klein{\min}} = \min \left\{ \frac{|d_h| \|\matAh^{-1}\|_{\infty}}
  {1-|e_h|\|\matAh^{-1}\|_{\infty}},  \frac{|d_h|\|\matAh^{-1}\|_{\infty}}
  {1-|d_h|\|\matAh^{-1}\|_{\infty}} \right\} =  \frac{|d_h|\|\matAh^{-1}\|_{\infty}}
  {1-|e_h|\|\matAh^{-1}\|_{\infty}} \leq 1.
\]
Similarly, from $|d_{H}|>|e_{H}|$ it follows that
\[
\eta_{H,\infty}^{\klein{\min}}
= \min \left\{\frac{|e_H|\|\matAH^{-1}\|_{\infty}}{1-|d_H| \|\matAH^{-1}\|_{\infty}},
\frac{|e_H|\|\matAH^{-1}\|_{\infty} }{1-|e_H|\|\matAH^{-1}\|_{\infty}}\right\}
= \frac{|e_H|\|\matAH^{-1}\|_{\infty} }{1-|e_H|\|\matAH^{-1}\|_{\infty}}
\leq \left|\frac{e_H}{d_H}\right|.
\]
Hence, an upper bound on $\rho_{12}$ and $\rho_{21}$ is given by
\[
\frac{\eta_{h,\infty}^{\klein{\min}} \|\matA^{-1}\matC\|_{\infty}}{1-\eta_{h,\infty}^{\klein{\min}}
\|\matA^{-1}\matB\|_{\infty}}\;
\frac{\eta_{H,\infty}^{\klein{\min}}\|\matA^{-1}\matB\|_{\infty}}{1-
\eta_{H,\infty}^{\klein{\min}}\|\matA^{-1}\matC\|_{\infty}}
\leq \eta_{h,\infty}^{\klein{\min}}\eta_{H,\infty}^{\klein{\min}} \leq \left|\frac{e_H}{d_H}\right| \equiv\rho.
\]
Finally,
\[
	\| T_{12} \|_{\infty} \leq
	\eta_{H,\klein{\infty}} \leq
	\left|\frac{e_H}{d_H}\right|,\qquad
	\| T_{21} \|_{\infty} \leq
	\eta_{h,\klein{\infty}}^{\klein{\min}} \leq 1,
\]
and substituting the values of $e_H$ and $d_H$ given in \eqref{eq:upwind}
yields the desired result.
\end{proof}

Note that the bound \eqref{eq:cor:conv-diff_1} does not depend on the
choice of $N$, i.e., the size of the mesh in the $x$-direction. Moreover, for
a fixed choice of $M$, and hence of $H_y$, the value of $\epsilon/(\epsilon+H_y)$
decreases with decreasing~$\epsilon$. Similar to the one-dimensional model
problem studied in~\cite{EchLieSzyTic18}, this indicates a faster convergence
of the multiplicative Schwarz method for smaller $\epsilon$, meaning larger
convection-dominance, which is confirmed in the following numerical example.

\begin{example}{\rm
We consider the model problem \eqref{eq:bvp2} with $\beta=0$, $f=0$,
and boundary conditions determined by the function
\begin{equation}\label{analytic_solution}
u(x,y)=(2x-1)\left(\frac{1-e^{(y-1)/\epsilon}}{1-e^{-1/\epsilon}}\right),
\end{equation}
which is a minor variation of the model problem in~\cite[Example 6.1.1]{ElmSilWat14}.
All computations in this example were done in MATLAB R2019a.
\Cref{fig:analytic_solution} shows $u(x,y)$ with $\epsilon=0.01$ on
$\Omega=(0,1)\times (0,1)$, and we clearly see the boundary layer close
to $y=1$, even for this relatively large value of $\epsilon$.

\begin{figure}[tbhp]
\hspace*{2em}
\begin{minipage}[t]{0.5\linewidth}
\includegraphics[width=0.83\linewidth]{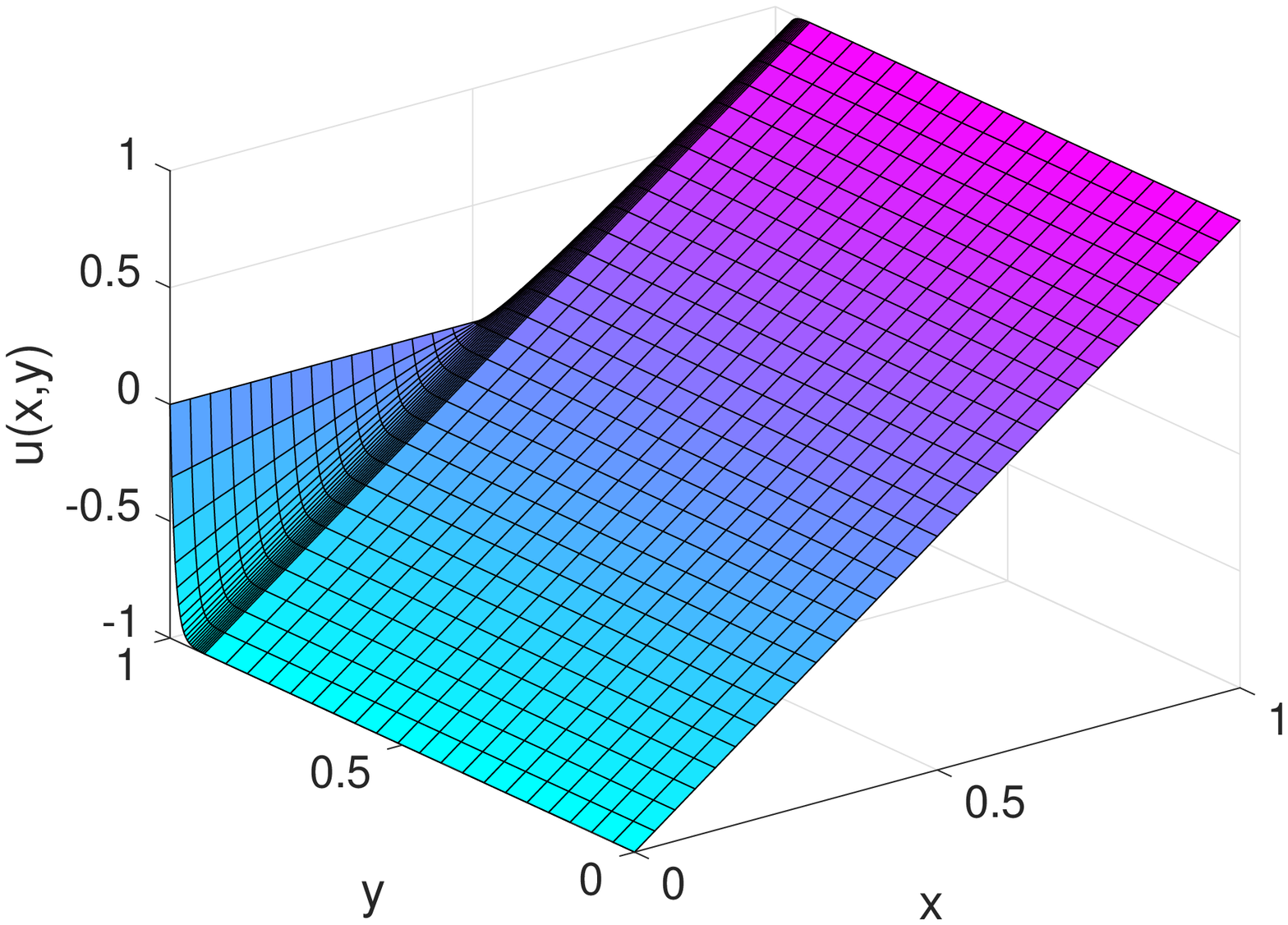}
\end{minipage}
\hspace*{-2.7em}
\begin{minipage}[t]{0.5\linewidth}
\includegraphics[width=0.83\linewidth]{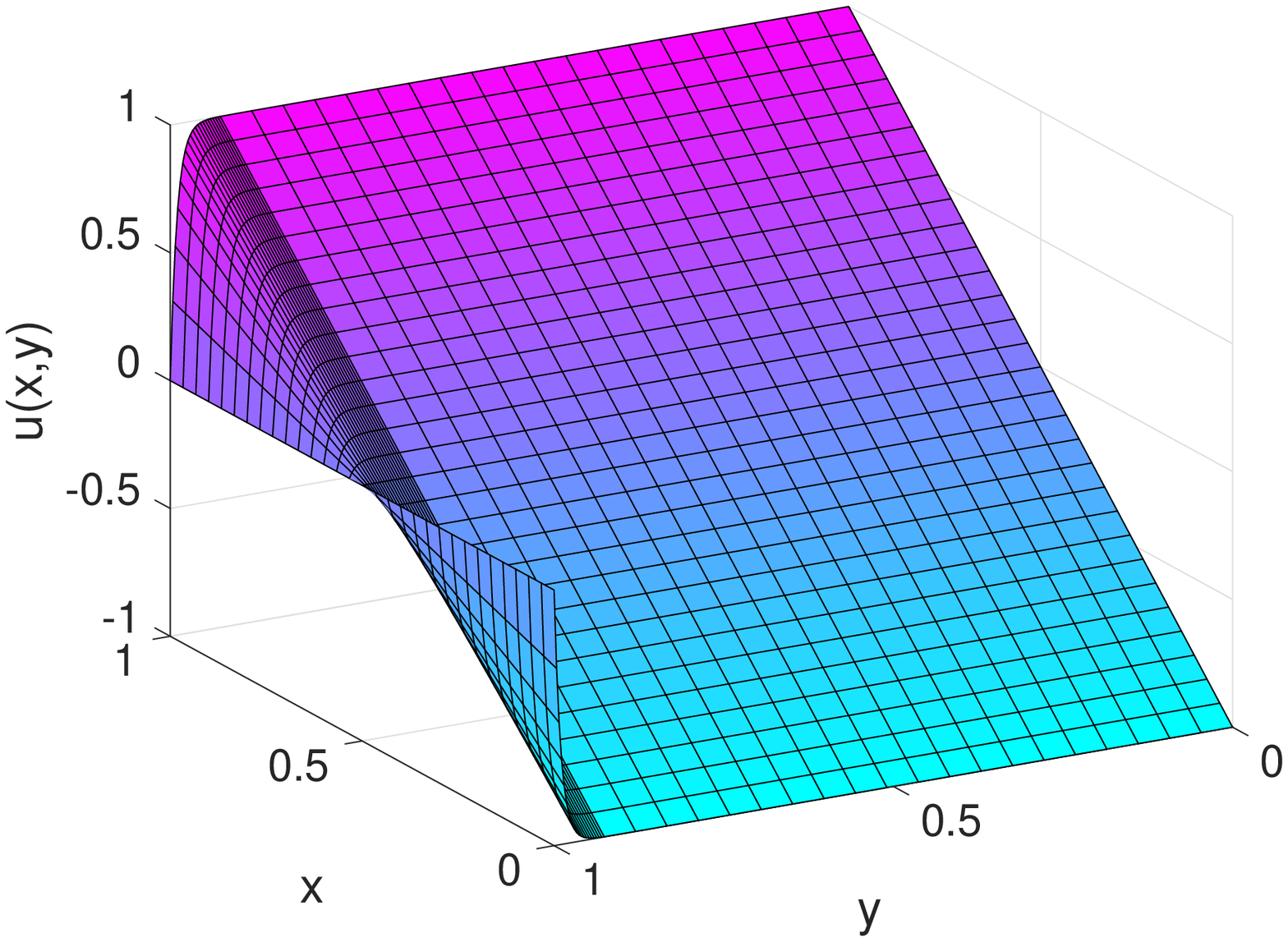}
\end{minipage}
\caption{The function \eqref{analytic_solution} with $\epsilon=0.01$ on
$\Omega=(0,1)\times (0,1)$ viewed from two different angles.}
\label{fig:analytic_solution}
\end{figure}

We discretize the problem using upwind differences on a Shishkin mesh
as described above with $N = 30$ and $M = 40$, and thus obtain a linear
algebraic system ${\cal A}x=b$ with ${\cal A}$ of size $1131\times 1131$.
We compute (an approximation to) the exact solution $x={\cal A}^{-1}b$
using MATLAB's backslash operator, and apply the multiplicative
Schwarz method with $x^{(0)}=0$. \Cref{fig:err1} shows the resulting
relative error norms
\[
\frac{\|x-x^{(k)}\|_\infty}{\|x\|_\infty},\;k=0,1,2,\ldots
\]
for the two iteration matrices $T_{12}$ and $T_{21}$ (solid lines)
and the corresponding upper bounds from \Cref{cor:conv-diff} for
$\epsilon=10^{-4}$ (left) and  $\epsilon=10^{-8}$ (right).
In all cases the bounds are very close to the actual error norms, and, as
indicated above, the multiplicative Schwarz method converges faster for
problems that are more convection-dominated.

\begin{figure}[tbhp]
\hspace*{-0.5em}
\begin{minipage}[t]{0.5\linewidth}
\includegraphics[width=0.95\linewidth]{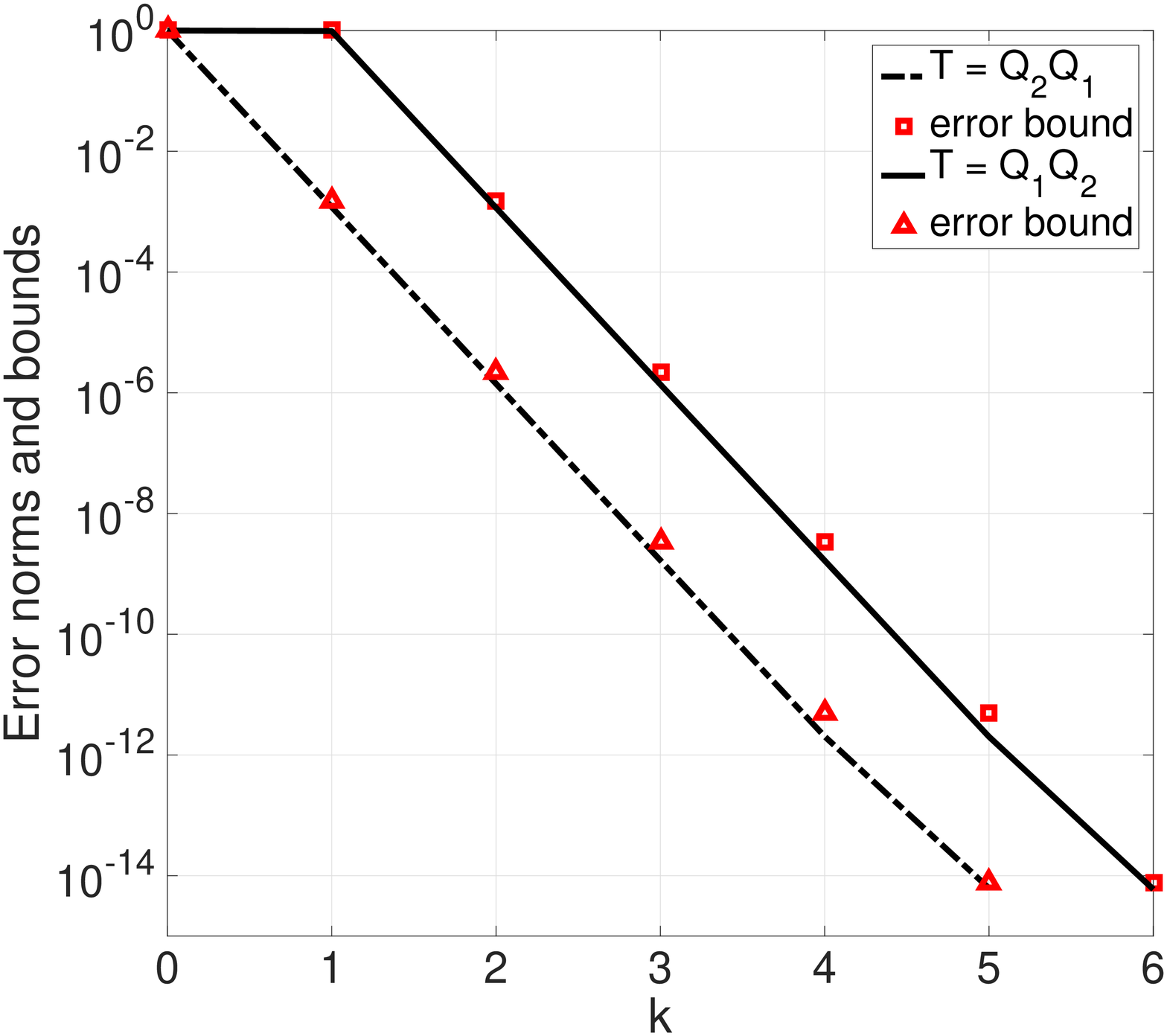}
\end{minipage}
\begin{minipage}[t]{0.5\linewidth}
\includegraphics[width=0.95\linewidth]{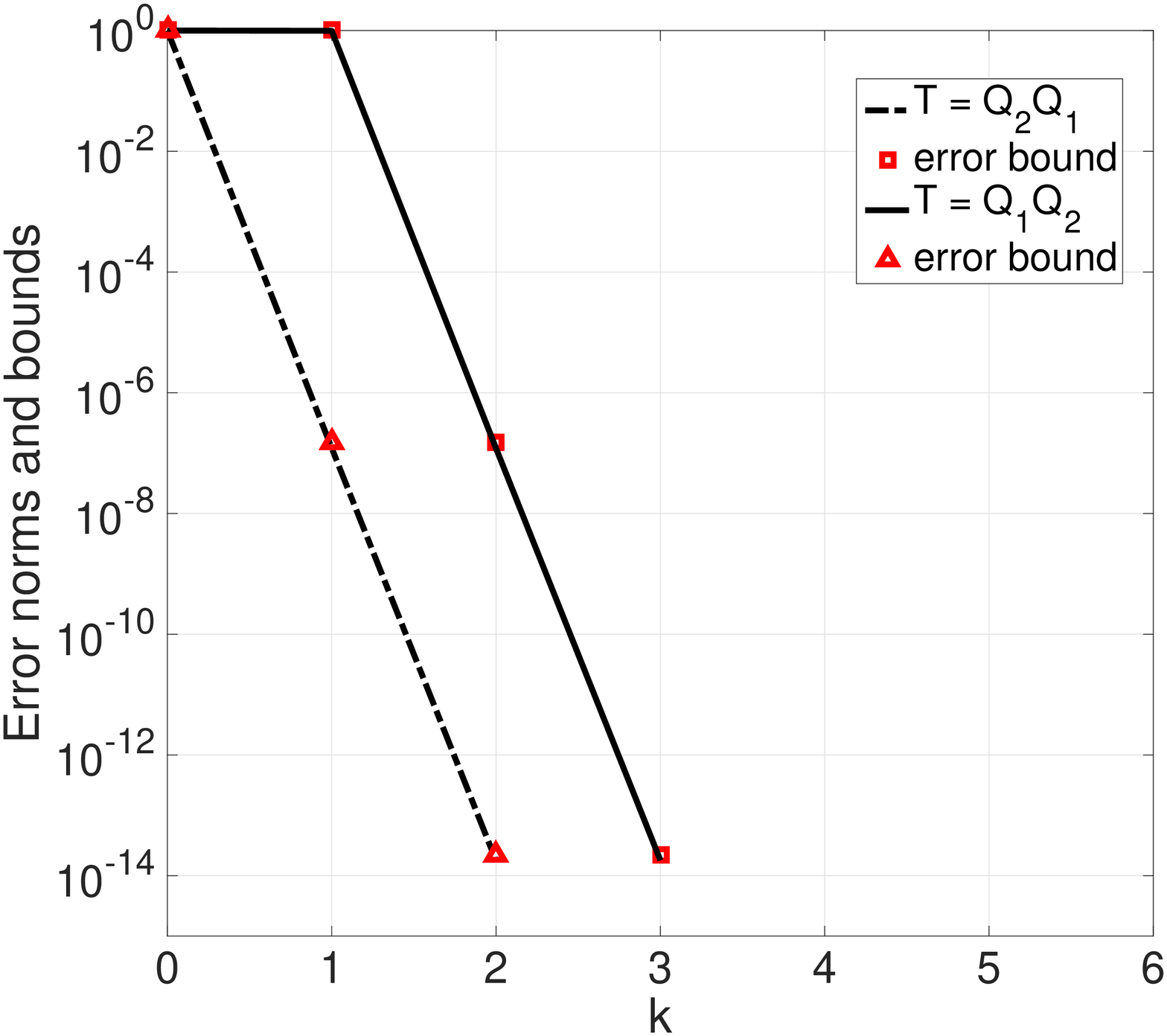}
\end{minipage}
\vspace*{-2em}
\caption{Convergence of multiplicative Schwarz and error bounds
for $\epsilon=10^{-4}$ (left) $\epsilon=10^{-8}$ (right).}\label{fig:err1}
\end{figure}

In \Cref{tab:2D:rho} we show, for different $N$, $M$ and $\epsilon$,
the values of the convergence factor $\rho_{12}$ as given in \eqref{eqn:rho} with
$\|\cdot\|=\|\cdot\|_{\infty}$, and the value of $\rho$ in
\eqref{eq:cor:conv-diff_2}, which represents an upper bound on $\rho_{12}$.
In all cases the two values are quite close to each other.

\begin{table}[tbhp]
\hspace*{1.3em}
\begin{minipage}[t]{0.5\linewidth}
\begin{tabular}{r|c c}
\cline{1-3}
\multicolumn{3}{c}{$N=20,M=20$, ${\cal
A}\in\mathbb{R}^{361\times361}$}\\\hline
\multicolumn{1}{ c|   }{$\epsilon$}& $\rho_{12}$ in~\eqref{eqn:rho} &
$\rho$ in~\eqref{eq:cor:conv-diff_2}\\
\multicolumn{1}{ c|  }{$10^{-8}$} & $7.5 \times 10^{-8}$ & $1.0 \times 10^{-7}$ \\
\multicolumn{1}{ c|  }{$10^{-6}$} & $7.5 \times 10^{-6}$ & $1.0 \times 10^{-5}$ \\
\multicolumn{1}{ c|  }{$10^{-4}$} & $7.5 \times 10^{-4}$ & $1.0 \times 10^{-3}$ \\
\multicolumn{1}{ c|  }{$10^{-2}$} & $7.0 \times 10^{-2}$ & $9.6 \times 10^{-2}$ \\ \hline
\multicolumn{3}{c}{$N=30, M=30$, ${\cal
A}\in\mathbb{R}^{841\times841}$}\\\hline
\multicolumn{1}{ c|  }{$10^{-8}$} & $1.2 \times 10^{-7}$ & $1.5 \times 10^{-7}$ \\
\multicolumn{1}{ c|  }{$10^{-6}$} & $1.2 \times 10^{-5}$ & $1.5 \times 10^{-5}$ \\
\multicolumn{1}{ c|  }{$10^{-4}$} & $1.2 \times 10^{-3}$ & $1.5 \times 10^{-3}$ \\
\multicolumn{1}{ c|  }{$10^{-2}$} & $1.1 \times 10^{-1}$ & $1.4 \times 10^{-1}$ \\ \hline
\multicolumn{3}{c}{$N=40,M=40$, ${\cal
A}\in\mathbb{R}^{1521\times1521}$}\\\hline
\multicolumn{1}{ c|  }{$10^{-8}$} & $1.7 \times 10^{-7}$ & $2.0 \times 10^{-7}$ \\
\multicolumn{1}{ c|  }{$10^{-6}$} & $1.7 \times 10^{-5}$ & $2.0 \times 10^{-5}$ \\
\multicolumn{1}{ c|  }{$10^{-4}$} & $1.7 \times 10^{-3}$ & $2.0 \times 10^{-3}$ \\
\multicolumn{1}{ c|  }{$10^{-2}$} & $1.4 \times 10^{-1}$ & $1.8 \times 10^{-1}$ \\\hline
\multicolumn{3}{c}{$N=50,M=50$, ${\cal
A}\in\mathbb{R}^{2401\times2401}$}\\\hline
\multicolumn{1}{ c|  }{$10^{-8}$} & $2.2 \times 10^{-7}$ & $2.5 \times 10^{-7}$ \\
\multicolumn{1}{ c|  }{$10^{-6}$} & $2.2 \times 10^{-5}$ & $2.5 \times 10^{-5}$ \\
\multicolumn{1}{ c|  }{$10^{-4}$} & $2.2 \times 10^{-3}$ & $2.5 \times 10^{-3}$ \\
\multicolumn{1}{ c|  }{$10^{-2}$} & $1.8 \times 10^{-1}$ & $2.1 \times 10^{-1}$ \\
\end{tabular}
\end{minipage}
\begin{minipage}[t]{0.5\linewidth}
\begin{tabular}{r|c c}
\cline{1-3}
\multicolumn{3}{c}{$N=20, M=30$, ${\cal
A}\in\mathbb{R}^{551\times551}$}\\\hline
\multicolumn{1}{ c|   }{$\epsilon$}& $\rho_{12}$ in~\eqref{eqn:rho} &
$\rho$ in~\eqref{eq:cor:conv-diff_2}\\
\multicolumn{1}{ c|  }{$10^{-8}$} & $1.2 \times 10^{-7}$ & $1.5 \times 10^{-7}$ \\
\multicolumn{1}{ c|  }{$10^{-6}$} & $1.2 \times 10^{-5}$ & $1.5 \times 10^{-5}$ \\
\multicolumn{1}{ c|  }{$10^{-4}$} & $1.2 \times 10^{-3}$ & $1.5 \times 10^{-3}$ \\
\multicolumn{1}{ c|  }{$10^{-2}$} & $1.1 \times 10^{-1}$ & $1.4 \times 10^{-1}$ \\ \hline
\multicolumn{3}{c}{$N=30,M=40$, ${\cal
A}\in\mathbb{R}^{1131\times1131}$}\\\hline
\multicolumn{1}{ c|  }{$10^{-8}$} & $1.7 \times 10^{-7}$ & $2.0 \times 10^{-7}$ \\
\multicolumn{1}{ c|  }{$10^{-6}$} & $1.7 \times 10^{-5}$ & $2.0 \times 10^{-5}$ \\
\multicolumn{1}{ c|  }{$10^{-4}$} & $1.7 \times 10^{-3}$ & $2.0 \times 10^{-3}$ \\
\multicolumn{1}{ c|  }{$10^{-2}$} & $1.4 \times 10^{-1}$ & $1.8 \times 10^{-1}$ \\\hline
\multicolumn{3}{c}{$N=40, M=50$, ${\cal
A}\in\mathbb{R}^{1911\times1911}$}\\\hline
\multicolumn{1}{ c|  }{$10^{-8}$} & $2.2 \times 10^{-7}$ & $2.5 \times 10^{-7}$ \\
\multicolumn{1}{ c|  }{$10^{-6}$} & $2.2 \times 10^{-5}$ & $2.5 \times 10^{-5}$ \\
\multicolumn{1}{ c|  }{$10^{-4}$} & $2.2 \times 10^{-3}$ & $2.5 \times 10^{-3}$ \\
\multicolumn{1}{ c|  }{$10^{-2}$} & $1.8 \times 10^{-1}$ & $2.1 \times 10^{-1}$ \\ \hline
\multicolumn{3}{c}{$N=50,M=60$, ${\cal
A}\in\mathbb{R}^{2891\times2891}$}\\\hline
\multicolumn{1}{ c|  }{$10^{-8}$} & $2.6 \times 10^{-7}$ & $3.0 \times 10^{-7}$ \\
\multicolumn{1}{ c|  }{$10^{-6}$} & $2.6 \times 10^{-5}$ & $3.0 \times 10^{-5}$ \\
\multicolumn{1}{ c|  }{$10^{-4}$} & $2.6 \times 10^{-3}$ & $3.0 \times 10^{-3}$ \\
\multicolumn{1}{ c|  }{$10^{-2}$} & $2.1 \times 10^{-1}$ & $2.5 \times 10^{-1}$ \\
\end{tabular}
\end{minipage}\\[1ex]
\caption{Values of $\rho_{12}$ computed using \eqref{eqn:rho} with
$\|\cdot\|=\|\cdot\|_{\infty}$ and $\rho$ in~\eqref{eq:cor:conv-diff_2}
for different values of $N$, $M$ and $\epsilon$.}
\label{tab:2D:rho}
\end{table}
}\end{example}

\section{Concluding discussion}\label{sec:conclusions}

Motivated by an analysis for a one-dimensional convection-diffusion model
problem in~\cite{EchLieSzyTic18}, we have studied the convergence of the
multiplicative Schwarz method for matrices with a special block structure.
After deriving a general expression for the convergence factor of the method,
we have focussed on block tridiagonal matrices, and applied recent results on
block diagonal dominance from~\cite{EchLieNab18} in order to obtain
quantitative error bounds that are valid from the first iteration.
In our analysis we did not use any of the usual assumptions on the matrices in
this context, such as symmetry, or the $M$- or $H$-matrix properties. We
illustrated our bounds numerically on a two-dimensional convection-diffusion
model problem that was discretized using a Shishkin mesh, and we found that
the bounds are very close to the actual error norms produced by the method.
We will now briefly discuss possible generalizations and alternative
applications of our approach.

\subsection*{Variable block sizes of ${\cal A}$ in \eqref{eq:blockmat}}
The results of \Cref{sec:schwarz}-\Cref{sec:bounds}, which are based on the
algebraic structure of the matrix \eqref{eq:blockmat} only, can be
generalized to matrices $\matAHhat$ and $\matAhhat$ having different
sizes. Such a generalization is straightforward if the sizes
of $\matAHhat$ and $\matAhhat$ are multiples of $N$, i.e.,
$\matAHhat\in \R^{Nm\times Nm}$ and $\matAhhat \in \R^{Ns\times Ns}$
with $s\neq m$, and more technical when this assumption is not satisfied.
We have chosen $\matAHhat$ and $\matAhhat$ of the same size in order to
reduce the technicalities in the analysis, and since our main application
in \Cref{sec:application} relies on the idea of the Shishkin mesh,
which is composed of piecewise equidistant meshes which have the same
number of unknowns.

\subsection*{General tridiagonal blocks in \eqref{eqn:tridiag}}
In \eqref{eqn:tridiag} we have assumed that $\matAHhat$ and $\matAhhat$
have a block Toeplitz structure. Apart from simplicity of notation, the
main motivation of this assumption was also the structure of the
discretized problem \eqref{eq:bvp2} with constant coefficients
studied in \Cref{sec:application}. Suppose that ${\cal A}$ has a
general block tridiagonal structure, i.e., that
\begin{equation}\label{eqn:tridiaggeneral}
\matAHhat=\mathrm{tridiag}(C_{\klein{H},i},A_{\klein{H},i},B_{\klein{H},i})
\quad\text{and}\quad
\matAhhat=\mathrm{tridiag}(C_{\klein{h},i},A_{\klein{h},i},B_{\klein{h},i}),
\end{equation}
for $i=1,\dots,m$. Such a matrix ${\cal A}$ can be obtained, for example,
by discretizing a boundary value problem like \eqref{eq:bvp2}, but with
nonconstant coefficients. An analysis of the multiplicative Schwarz method
for such a matrix ${\cal A}$ following the approach in this paper is
still possible, since the results from~\cite{EchLieNab18} on block
diagonal dominance (see also Appendix A) are formulated for general
block tridiagonal matrices. A generalization of \Cref{thm:main2}
to ${\cal A}$ with blocks \eqref{eqn:tridiaggeneral} would require that
the conditions \eqref{eqn:dominant} hold in every block row, and then
analogously to \eqref{eqn:etah}--\eqref{eqn:etaH}, every block row in
$\matAHhat$ or $\matAhhat$ would give a parameter $\eta_{H,i}$
or $\eta_{h,i}$, respectively.

\subsection*{Problems with two boundary layers}
\begin{figure}
\hspace*{-2em}
\begin{minipage}[b]{0.5\textwidth}
\begin{center}
\begin{tikzpicture}[scale = 0.4]
\draw (0,0) rectangle (10, 10);
\node[left] at (0,10) {$1$};
\draw (0,8) node[left] {$1-\tau_y$} -- (10,8);
\draw (8,0) node[below] {$1-\tau_x$} -- (8,10);
\node[left] at (0,0) {$0$};
\node[below] at (0,0) {$0$};
\node[below] at (10,0) {$1$};
\draw[black,fill=black] (0,8) circle (0.5ex);
\node at (4.2,4) {$\Omega_{11}$};
\node at (4.2,9) {$\Omega_{21}$};
\node at (9.1,9) {$\Omega_{22}$};
\node at (9.1,4) {$\Omega_{12}$};
\end{tikzpicture}
\end{center}
\end{minipage}
\hspace*{1em}
\begin{minipage}[b]{0.5\textwidth}
\begin{center}
\begin{tikzpicture}[scale = 0.4]
\draw (0,0) rectangle (10, 10);

\foreach \point in {2,4,6,8,8.5,9,9.5}
\draw (0,\point) -- (10,\point);
\foreach \point in {2,4,6,8,8.5,9,9.5}
\draw (\point,0) -- (\point,10);

\node[left] at (0,10) {$y_{\klein{M}}$};
\node[left] at (0,8)  {$y_{\klein{\frac{M}{2}}}$};
\node[left] at (0,0)  {$y_0$};

\node[below] at (10.4,0) {$x_{\klein{M}}$};
\node[below] at (8.4,0)  {$x_{\klein{\frac{M}{2}}}$};
\node[below] at (0.2,0)  {$x_0$};

\foreach \point in {2,4,6,8}
\draw (\point,0) -- (\point,10);

\draw[black,fill=black] (0,8) circle (0.5ex);
\draw (5,-1.25) node[above] {$H_x$};
\draw (0,5) node[left] {$H_y$};
\draw (10,9) node[right] {$h_y$};
\draw (9,10) node[above] {$h_x$};
\end{tikzpicture}
\end{center}
\end{minipage}
\caption{Division of the domain and Shishkin mesh for the problem
\eqref{eq:bvp3} with two boundary layers.}
\label{fig:2Ddomain3}
\end{figure}
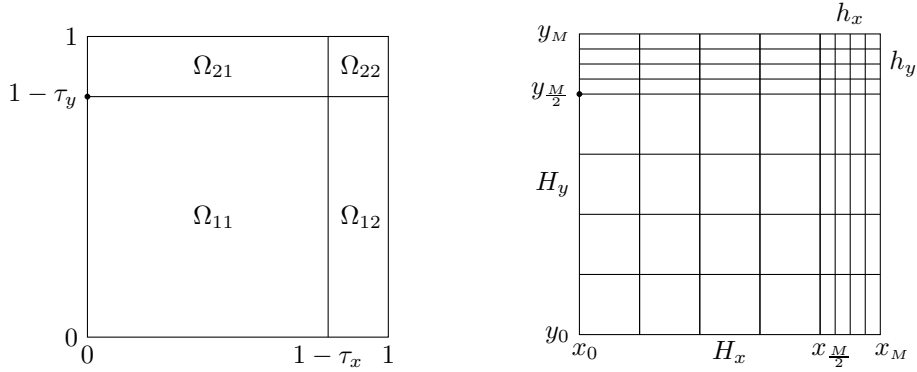

More practical two-dimensional problems arise when considering convection
in both directions, e.g.,
\begin{equation}\label{eq:bvp3}
-\epsilon  \Delta u + u_x + u_y +\beta u = f  \text{ in } \Omega=(0,1)\times(0,1),
\;0<\epsilon\ll 1,\;\beta\geq 0,
\;\; u=0\,\text{  on  }\, \partial\Omega.
\end{equation}
In this case the solution has two boundary layers at $x=1$ and $y=1$;
for more details, see, e.g.,~\cite{LinSty01,MilOriShi96,RooStyTob08}.
One can again use the appropriate Shishkin mesh to resolve the boundary
layers (see \Cref{fig:2Ddomain3}) and discretize the problem using
standard upwind finite differences and the lexicographical line ordering
of the unknowns (see, e.g.,~\cite{LinSty01}), which yields a matrix
$\mathcal A$ with a block tridiagonal structure. The individual
blocks are not Toeplitz any more, but their more complicated
``Toeplitz-like'' structure can be analyzed using the chosen
discretization scheme.

To solve the resulting algebraic system we can again apply the multiplicative
Schwarz method. Now we have four regions $\Omega_{ij}$ (see \Cref{fig:2Ddomain3}), and we can choose restriction operators corresponding to these regions.
Our numerical experiments predict that the multiplicative Schwarz
method converges for any ordering of subdomains, but some orderings
lead to a faster convergence (e.g., $\Omega_{11}$, $\Omega_{12}$, $\Omega_{21}$, $\Omega_{22}$) than the others
(e.g., $\Omega_{12}$, $\Omega_{22}$, $\Omega_{21}$, $\Omega_{11}$).
In all cases the iteration matrices seem to have a low numerical rank
(close to $N$). The algebraic structure of the iteration matrices
is now more complicated, but we believe that it is still analyzable
along the lines of \Cref{sec:structure}. A full analysis is, however,
beyond the scope of this paper.

\subsection*{Multiplicative Schwarz preconditioned GMRES}
It can be seen from \eqref{eq:schwarz} that the fixed point of the multiplicative
Schwarz iteration, and hence the solution of \eqref{eq:linsys}, satisfies
\begin{equation}\label{eq:precond}
(I-T_{ij})x=v.
\end{equation}
As discussed also in~\cite[Section~6]{EchLieSzyTic18}, if we apply the Krylov subspace method GMRES~\cite{SaaSch86} to \eqref{eq:precond}, then the multiplicative
Schwarz method can be seen as a preconditioner; see, e.g., \cite{KahKamPhi07}. For the matrices
studied in this paper we have ${\rm rank}(T_{ij})\leq N$ (see \Cref{lem:powers}),
and therefore
\[
\dim \left(\mathcal{K}_k(I-T_{ij},r_0)\right) \leq N+1,
\]
for any initial residual $r_0$. Consequently, the multiplicative Schwarz
GMRES will converge to the solution in at most $N+1$ steps (in exact
arithmetic), even when the multiplicative Schwarz method itself converges
slowly or diverges. Note that $T_{ij}$ will also have low rank for
discretized boundary value problem like \eqref{eq:bvp2} with
nonconstant coefficients, and possibly also in more complicated problems
with several boundary layers.

\subsection*{Additive Schwarz method} In this paper we obtained results for the  multiplicative Schwarz method solving linear algebraic systems with $\mathcal{A}$ given by \eqref{eq:blockmat}. A natural question is whether analogous results
can be obtained when considering the additive Schwarz method. In that method
we use the iteration scheme
\[
x^{(k+1)}=Tx^{(k)}+(P_{1}+P_{2})b,\quad T\equiv I-(P_{1}+P_{2});
\]
see, e.g., \cite{BenFroNabSzy01,NabSzy06}. Using the structure of projection
matrices $P_1$ and $P_2$ described in \Cref{sec:structure} we obtain
\begin{eqnarray*}
T = -\left[\begin{array}{ccccc}
0_{N(m-1)} &  &  & P_{1:m-1}^{(1)}\\
  &  &  & P_{m}^{(1)}\\
  & \Pi^{(2)} & I_{N} & \Pi^{(1)}\\
  & P_{1}^{(2)} &  & \\
  & P_{2:m}^{(2)} &  &  & 0_{N(m-1)}
\end{array}\right].
\end{eqnarray*}
Obviously, the spectral radius of $T$ is larger than or equal to one, so that
the additive Schwarz method is not convergent in general.
Moreover, the matrix $I-T$ in the preconditioned system $(I-T)x=v$ is singular,
so that the additive Schwarz method cannot be used as a preconditioner
for GMRES in general. The same result can be obtained also for the damped
additive Schwarz method considered in \cite{NabSzy06}.

\pagebreak
\appendix
\label{sec:Appendix}

\section{Column block diagonal dominance of matrices}\label{app:colBDD}\\
Analogous to \cite[Definition~2.1]{EchLieNab18} we can define column
block diagonal dominance of matrices as follows.

\begin{definition}\label{def:coldom}
A matrix of the form
\[
A=[A_{ij}]\quad
\mbox{with blocks $A_{ij}\in\mathbb{C}^{m\times m}$ for $i,j=1,\dots,n$}
\]
is called \emph{column block diagonally dominant} (with respect
to the matrix norm $\|\cdot\|$), when the diagonal blocks $A_{jj}$
are nonsingular, and
\begin{equation}\label{eq:blocdiagdom1}
\sum_{\atopfrac{i=1}{i\neq j}}^{n} \|A_{ij}A_{jj}^{-1}\| \leq 1,
\quad \text{for $j=1,\dots,n$}.
\end{equation}
If strict inequality holds in \eqref{eq:blocdiagdom1} then $A$ is called
\emph{column block strictly diagonally dominant} (with respect to the matrix
norm $\|\cdot\|$).
\end{definition}

Restricting our attention to block \emph{tridiagonal} matrices
$A={\rm tridiag}(C_i,A_i,B_i)$ as in \cite[Equation (2.4)]{EchLieNab18},
and following the notation of that paper, we define
\[
\tilde{\tau}_i  = \frac{\|C_iA_i^{-1}\|}{1-\|B_{i-1}A_i^{-1}\|},
\quad
\tilde{\omega}_i = \frac{\|B_{i-1}A_i^{-1}\|}{1-\|C_{i}A_i^{-1}\|},
\quad \text{for $i=1,\ldots,n$},
\]
where $B_0=C_n=0$. The column block diagonal dominance of $A$ then
implies that $0\leq \tilde{\tau}_i\leq 1$ and $0\leq \tilde{\omega}_i \leq 1$.
The same approach as in \cite[Section~2]{EchLieNab18} now yields the
the following result, which is the ``column version'' of~\cite[Theorem~2.6]{EchLieNab18}.
\begin{theorem}\label{thm:coldom}
Let $A={\rm tridiag}(C_i,A_i,B_i)$ be column block diagonally dominant, and
suppose that the blocks $B_i$, $C_i$ for $i=1,\ldots,n-1$ are nonsingular. If
in addition
\[
\|C_1A_1^{-1}\|<1 \quad \mbox{and} \quad \|B_{n-1}A_n^{-1}\|<1,
\]
then $A^{-1}=[Z_{ij}]$ with
\begin{eqnarray*}
\|Z_{ij}\| & \leq & \|Z_{ii}\| \prod_{k=i+1}^{j}\tilde{\omega}_k,
\quad\text{for all $i<j$}, \\ 
\|Z_{ij}\| & \leq & \|Z_{ii}\| \;\,\prod_{k=j}^{i-1}\tilde{\tau}_k,
\qquad\text{for all $i>j$}. 
\end{eqnarray*}
Moreover, for $i=1,\dots,n$,
\[
\frac{\|I\|}
{\|A_i\|+\tilde{\tau}_{i-1}\|B_{i-1}\|+\tilde{\omega}_{i+1}\|C_{i}\|}
\leq
\|Z_{ii}\|
\leq
\frac{\|I\|}
{\|A_i^{-1}\|^{-1}-\tilde{\tau}_{i-1}\|B_{i-1}\|-\tilde{\omega}_{i+1}\|C_i\|},
\]
provided that the denominator of the upper bound is larger than zero,
and where we set $\tilde{\tau}_0=\tilde{\omega}_{n+1}=0$.
\end{theorem}

\bibliographystyle{siam}
\bibliography{EchLieTic20}

\begin{thebibliography}{10}

\bibitem{BenFroNabSzy01}
{\sc M.~Benzi, A.~Frommer, R.~Nabben, and D.~B. Szyld}, {\em Algebraic theory
  of multiplicative {S}chwarz methods}, Numer. Math., 89 (2001), pp.~605--639.

\bibitem{BruPedSzy04}
{\sc R.~Bru, F.~Pedroche, and D.~B. Szyld}, {\em Overlapping additive and
  multiplicative {S}chwarz iterations for {$H$}-matrices}, Linear Algebra
  Appl., 393 (2004), pp.~91--105.

\bibitem{EchLieNab18}
{\sc C.~Echeverr\'{i}a, J.~Liesen, and R.~Nabben}, {\em Block diagonal
  dominance of matrices revisited: bounds for the norms of inverses and
  eigenvalue inclusion sets}, Linear Algebra Appl., 553 (2018), pp.~365--383.

\bibitem{EchLieSzyTic18}
{\sc C.~Echeverr{\'i}a, J.~Liesen, D.~B. Szyld, and P.~Tich{\'y}}, {\em
  Convergence of the multiplicative {S}chwarz method for singularly perturbed
  convection-diffusion problems discretized on a {S}hishkin mesh}, Electron.
  Trans. Numer. Anal., 48 (2018), pp.~40--62.

\bibitem{ElmSilWat14}
{\sc H.~C. Elman, D.~J. Silvester, and A.~J. Wathen}, {\em Finite elements and
  fast iterative solvers: with applications in incompressible fluid dynamics},
  Oxford University Press, Oxford, second~ed., 2014.

\bibitem{FroNabSzy08}
{\sc A.~Frommer, R.~Nabben, and D.~B. Szyld}, {\em Convergence of stationary
  iterative methods for {H}ermitian semidefinite linear systems and
  applications to {S}chwarz methods}, SIAM J. Matrix Anal. Appl., 30 (2008),
  pp.~925--938.

\bibitem{FroSzy14}
{\sc A.~Frommer and D.~B. Szyld}, {\em On necessary conditions for convergence
  of stationary iterative methods for {H}ermitian semidefinite linear systems},
  Linear Algebra Appl., 453 (2014), pp.~192--201.

\bibitem{HornJohn12}
{\sc R.~A. Horn and C.~R. Johnson}, {\em Matrix Analysis}, Cambridge University
  Press, Cambridge, second~ed., 2013.

\bibitem{KahKamPhi07}
{\sc G.~A.~A. Kahou, E.~Kamgnia, and B.~Philippe}, {\em An explicit formulation
  of the multiplicative {S}chwarz preconditioner}, Appl. Numer. Math., 57
  (2007), pp.~1197--1213.

\bibitem{KopOri10}
{\sc N.~Kopteva and E.~O'Riordan}, {\em Shishkin meshes in the numerical
  solution of singularly perturbed differential equations}, Int. J. Numer.
  Anal. Model., 7 (2010), pp.~393--415.

\bibitem{KraParSte83}
{\sc D.~Kratzer, S.~V. Parter, and M.~Steuerwalt}, {\em Block splittings for
  the conjugate gradient method}, Comput. \& Fluids, 11 (1983), pp.~255--279.

\bibitem{LinSty01}
{\sc T.~Lin\ss and M.~Stynes}, {\em Numerical methods on {S}hishkin meshes for
  linear convection-diffusion problems}, Comput. Methods Appl. Mech. Engrg.,
  190 (2001), pp.~3527--3542.

\bibitem{MilOriShi96}
{\sc J.~J.~H. Miller, E.~O'Riordan, and G.~I. Shishkin}, {\em Fitted Numerical
  Methods for Singular Perturbation Problems: Error Estimates in the Maximum
  Norm for Linear Problems in One and Two Dimensions}, World Scientific
  Publishing Co. Pte. Ltd., Hackensack, NJ, revised~ed., 2012.

\bibitem{NabSzy06}
{\sc R.~Nabben and D.~B. Szyld}, {\em Schwarz iterations for symmetric positive
  semidefinite problems}, SIAM J. Matrix Anal. Appl., 29 (2006/07),
  pp.~98--116.

\bibitem{RooStyTob08}
{\sc H.-G. Roos, M.~Stynes, and L.~Tobiska}, {\em Robust numerical methods for
  singularly perturbed differential equations}, vol.~24 of Springer Series in
  Computational Mathematics, Springer-Verlag, Berlin, second~ed., 2008.

\bibitem{SaaSch86}
{\sc Y.~Saad and M.~H. Schultz}, {\em G{MRES}: a generalized minimal residual
  algorithm for solving nonsymmetric linear systems}, SIAM J. Sci. Statist.
  Comput., 7 (1986), pp.~856--869.

\bibitem{Sty05}
{\sc M.~Stynes}, {\em Steady-state convection-diffusion problems}, Acta Numer.,
  14 (2005), pp.~445--508.

\end{thebibliography}
\end{document}